\documentclass[12pt]{article}

\usepackage{amsmath, amssymb, booktabs, color, algorithm, algorithmic}
\usepackage{amsthm}
\usepackage{bm}
\usepackage{caption}
\usepackage{subcaption}
\usepackage{here}
\usepackage{comment}
\usepackage[dvipdfmx]{graphicx}

\newtheorem{thm}{Theorem}[section]
\newtheorem{prop}{Proposition}[section]
\newtheorem{lem}{Lemma}[section]
\newtheorem{cor}{Corollary}[section]

\newtheorem{prob}{Problem}[section]
\newtheorem{algo}{Algorithm}[section]
\newtheorem{assum}{Assumption}[section]

\newcommand{\argmin}{\operatornamewithlimits{argmin}}

\numberwithin{equation}{section}

\begin{document}
\makeatletter

\begin{center}
\large{\bf Proximal Point Algorithms for Nonsmooth Convex Optimization with Fixed Point Constraints}\\
\small{This work was supported by the Japan Society for the Promotion of Science through a Grant-in-Aid for
Scientific Research (C) (15K04763).}
\end{center}\vspace{3mm}

\begin{center}
\textsc{Hideaki Iiduka}\\
Department of Computer Science, 
Meiji University,
1-1-1 Higashimita, Tama-ku, Kawasaki-shi, Kanagawa 214-8571 Japan. 
(iiduka@cs.meiji.ac.jp)
\end{center}

\vspace{2mm}

\footnotesize{
\noindent\begin{minipage}{14cm}
{\bf Abstract:}
The problem of minimizing the sum of nonsmooth, convex objective functions defined on a real Hilbert space over the intersection of fixed point sets of nonexpansive mappings, onto which the projections cannot be efficiently computed, is considered. The use of proximal point algorithms that use the proximity operators of the objective functions and incremental optimization techniques is proposed for solving the problem. With the focus on fixed point approximation techniques, two algorithms are devised for solving the problem. One blends an incremental subgradient method, which is a useful algorithm for nonsmooth convex optimization, with a Halpern-type fixed point iteration algorithm. The other is based on an incremental subgradient method and the Krasnosel'ski\u\i-Mann fixed point algorithm. It is shown that any weak sequential cluster point of the sequence generated by the Halpern-type algorithm belongs to the solution set of the problem and that there exists a weak sequential cluster point of the sequence generated by the Krasnosel'ski\u\i-Mann-type algorithm, which also belongs to the solution set. Numerical comparisons of the two proposed algorithms with existing subgradient methods for concrete nonsmooth convex optimization show that the proposed algorithms achieve faster convergence.
\end{minipage}
 \\[5mm]

}

\hbox to14cm{\hrulefill}\par


\section{Introduction}\label{sec:1}
Convex optimization theory is a powerful tool for solving many practical problems in operational research (see, e.g., 
\cite{ben2008,plate2010} and references therein).
In particular, it has been widely used to solve practical convex minimization problems over complicated constraints, e.g., convex optimization problems with a {\em fixed point constraint} \cite{com,iiduka_siopt2013,iiduka_hishinuma_siopt2014,mainge2014,yamada,yao2011} and with a {\em variational inequality constraint} \cite{facc2014,mainge2010_1}. 

Consider the following convex optimization problem: given a convex objective function $f \colon H \to \mathbb{R}$ and a nonexpansive mapping $T \colon H \to H$,
\begin{align}\label{problem:min}
\text{minimize } f(x) \text{ subject to } x \in \mathrm{Fix}(T),
\end{align}
where $H$ is a real Hilbert space and $\mathrm{Fix}(T)$ stands for the fixed point set of $T$. 
Problem \eqref{problem:min} enables consideration of optimization problems with complicated constraint sets 
\cite[section I]{com1999}, \cite[subsection 3.2]{neto2009}, \cite[section 4]{yamada} onto which metric projections cannot be easily calculated. Several algorithms (e.g., \cite{com,iiduka_siopt2013,iiduka_hishinuma_siopt2014,yamada}) have been proposed for solving problem \eqref{problem:min} when $f$ is smooth and convex, which includes practical problems such as signal recovery \cite{com}, beamforming \cite{slav}, and network resource allocation \cite{iiduka_siopt2013,iiduka_hishinuma_siopt2014}.

Here, problem \eqref{problem:min} is considered for when $f$ is convex but {\em not always smooth}. One objective is to devise optimization algorithms for nonsmooth convex optimization problem \eqref{problem:min}, which cannot be solved using conventional algorithms for smooth convex optimization \cite{com,iiduka_siopt2013,iiduka_hishinuma_siopt2014,yamada}.
There are significant problems with problem \eqref{problem:min} when $f$ is a general nonsmooth convex function (e.g., the $L^1$-norm). They include the problem of minimizing the total variation of a signal over a convex set, Tykhonov-like problems with $L^1$-norms \cite[I. Introduction]{comb2007}, the classifier ensemble problem with sparsity and diversity learning \cite[subsection 2.2.3]{yin2014_1}, \cite[subsection 3.2.4]{yin2014_2}, which is expressed as $L^1$-norm minimization, and the minimal antenna-subset selection problem \cite[subsection 17.4]{yamada2011}. Another objective is to solve problem \eqref{problem:min} including the above real-world problems by using {\em incremental} optimization techniques. If the explicit forms of $f$ and $T$ in problem \eqref{problem:min} are unknowable, algorithms making the best use of their mapping information cannot be applied to the problem. To enable us to consider such a case, a networked system with a finite number of users is assumed, and each user $i$ is assumed to have its own private convex, nonsmooth objective function $f^{(i)}$ and nonexpansive mapping $T^{(i)}$. The main objective is to devise optimization algorithms that enable each user to find an optimal solution to problem \eqref{problem:min} with
\begin{align}\label{FT}
f := \sum_{i=1}^I f^{(i)} \text{ and } \mathrm{Fix}(T) := \bigcap_{i=1}^I \mathrm{Fix} \left(T^{(i)} \right),
\end{align}
where $I\in \mathbb{N}$ is the number of users, without using the private information of other users.

There have been many reports on incremental and parallel optimization algorithms. Parallel proximal algorithms \cite[Proposition 27.8]{b-c}, which use the proximity operators of nonsmooth, convex functions, are useful for minimizing the sum of nonsmooth, convex functions over the whole space. Incremental subgradient methods \cite{neto2009,nedic2001,solo1998} and projected multi-agent algorithms \cite{nedic2009,nedic2009_1,nedic} can minimize the sum of nonsmooth, convex functions for certain constraint sets by using the subgradients of the nonsmooth, convex functions instead of the proximity operators. The incremental subgradient algorithm \cite{neto2009} and the asynchronous proximal algorithm \cite{pesq2015} can work on nonsmooth convex optimization over sublevel sets of convex functions onto which the projections cannot be easily calculated. The incremental and parallel gradient algorithms \cite{iiduka_siopt2013,iiduka_hishinuma_siopt2014} can work on smooth convex optimization over fixed point sets of nonexpansive mappings. The incremental and parallel algorithms \cite{iiduka_fpta2015,iiduka_mp2015,iiduka_2014} use the subgradients of nonsmooth convex functions and can optimize the sum of the nonsmooth convex functions over fixed point sets of nonexpansive mappings. To the best of our knowledge, there have been no reports on {\em incremental proximal point algorithms for nonsmooth convex optimization with fixed point constraints}.

Ideas from three useful types of algorithms, (I) proximal point algorithms, (II) incremental subgradient algorithms, and (III) fixed point algorithms, are used to achieve the main objective. 

(I) 
The well-known {\em proximal point algorithms} 
(see, e.g.,  \cite[Chapter 27]{b-c}, \cite{lions,martinet1970,rock1976} and references therein)
for nonsmooth convex optimization use the {\em proximity operators} \cite[Definition 12.23]{b-c}, \cite{moreau1962} of convex functions. Here, it is assumed that user $i$ can use the proximity operator of $f^{(i)}$, which is defined for all $x\in H$ by 
\begin{align*}
\mathrm{Prox}_{f^{(i)}} (x) \in \argmin_{y\in H} \left[ f^{(i)}(y) + \frac{1}{2} \left\|x-y\right\|^2  \right].
\end{align*}

(II) {\em Incremental subgradient algorithms} \cite{neto2009,nedic2001,solo1998} are useful
algorithms for nonsmooth convex optimization. An iteration $n$ of the algorithm is defined as follows: given $x_n^{(0)} \in H$,
\begin{align}\label{incremental}
\begin{split}
&x_n^{(i)} := x_n^{(i)} \left( x_n^{(i-1)}, f^{(i)}, T^{(i)} \right) \text{ } (i=1,2,\ldots,I),\\
&x_{n+1} := x_n^{(I)} =: x_{n+1}^{(0)}.
\end{split}
\end{align}
Under the assumption that user $i$ can communicate with neighbor user $(i-1)$, user $i$ can implement algorithm \eqref{incremental} by using only its own private mappings $f^{(i)},T^{(i)}$ and information $x_n^{(i-1)}$ transmitted from the neighbor user.

(III) There are many {\em fixed point algorithms} \cite{berinde} for solving fixed point problems. Here, the focus is on using the {\em Halpern fixed point algorithm} \cite{halpern,wit} and the {\em Krasnosel'ski\u\i-Mann fixed point algorithm} \cite{kra,mann} to search for a fixed point of a nonexpansive mapping $T$. The former is defined as follows: for each $n\in\mathbb{N}$, $x_{n+1} := \alpha_n x_0 + (1-\alpha_n ) T (x_n)$. The latter is defined as $x_{n+1} := \alpha_n x_n + (1-\alpha_n ) T (x_n)$, where $x_0\in H$ and $(\alpha_n)_{n\in\mathbb{N}} \subset [0,1]$. When user $i$ has $x^{(i)}\in H$, $f^{(i)}$, and $T^{(i)}$ and information $x_n^{(i-1)}$ transmitted from user $(i-1)$, user $i$ can compute
\begin{align}\label{H}
x_n^{(i)} := \alpha_n x^{(i)} + (1-\alpha_n) T^{(i)} \left(y_n^{(i)} \left(x_n^{(i-1)}, f^{(i)}\right) \right),
\end{align}
which is based on the Halpern fixed point algorithm, or 
\begin{align}\label{KM}
x_n^{(i)} := \alpha_n x_n^{(i-1)} + (1-\alpha_n) T^{(i)} \left(y_n^{(i)} \left(x_n^{(i-1)}, f^{(i)}\right) \right),
\end{align}
which is based on the Krasnosel'ski\u\i-Mann fixed point algorithm, where $y_n^{(i)}$ is a point depending on only $x_n^{(i-1)}$ and $f^{(i)}$. From (I), $y_n^{(i)}$ can be defined using the value of the proximity operator of $f^{(i)}$ at $x_n^{(i-1)}$; i.e.,
\begin{align}\label{prox}
y_n^{(i)} := \mathrm{Prox}_{f^{(i)}} \left(x_n^{(i-1)} \right).
\end{align}

Two incremental proximal point algorithms are proposed for solving problem \eqref{problem:min} with $f$ and $T$ defined by \eqref{FT}. One is based on the proximal point algorithm \eqref{prox}, the incremental subgradient method \eqref{incremental}, and the Halpern fixed point algorithm \eqref{H}. The other uses the ideas of the proximal point algorithm \eqref{prox}, the incremental subgradient method \eqref{incremental}, and the Krasnosel'ski\u\i-Mann fixed point algorithm \eqref{KM}.
 
Here, let us explicitly compare the two proposed algorithms with the existing algorithms \cite{iiduka_fpta2015,iiduka_mp2015,iiduka_2014}. The proposed and existing algorithms can be applied to problem \eqref{problem:min} with \eqref{FT}. The convergence analyses in \cite{iiduka_fpta2015,iiduka_mp2015,iiduka_2014} showed that {\em there exists a weak sequential cluster point} of the sequence generated by one of the existing algorithms that belongs to the solution set of problem \eqref{problem:min} with \eqref{FT}. However, these results are not strong enough. This is because knowing the existence of one optimal cluster point cannot help users to identify an optimal solution when multiple cluster points are observed. In contrast to the results in \cite{iiduka_fpta2015,iiduka_mp2015,iiduka_2014}, one of the proposed algorithms (Algorithm \ref{algorithm:2}) satisfies a gratifying convergence property such that {\em any weak sequential cluster point} of the sequence generated by the proposed algorithm belongs to the solution set of problem \eqref{problem:min} with \eqref{FT} under certain assumptions (Theorem \ref{theorem:2}). This result is attributed to the framework of the algorithm being based on the Halpern fixed point algorithm \eqref{H} as compared with the existing algorithms, which are based on the Krasnosel'ski\u\i-Mann fixed point algorithm \eqref{KM}. Since the other proposed algorithm (Algorithm \ref{algorithm:KM}) is based on the Krasnosel'ski\u\i-Mann fixed point algorithm \eqref{KM}, it is not guaranteed that any weak sequential cluster point of the sequence in Algorithm \ref{algorithm:KM} belongs to the solution set (Theorem \ref{theorem:KM}). However, Algorithm \ref{algorithm:KM} can work when the step sizes are constant, which are the most tractable choice of step size sequences, in contrast to Algorithm \ref{algorithm:2}, which uses diminishing step size sequences.

One contribution of this paper is analysis of the proposed algorithms' convergence. It is shown that, under certain assumptions, any weak sequential cluster point of the sequence generated by the Halpern-type algorithm belongs to the solution set of the problem and that there exists a weak sequential cluster point of the sequence generated by the Krasnosel'ski\u\i-Mann-type algorithm, which also belongs to the solution set. Another contribution of this paper is provision of examples showing that the proposed algorithms perform better than subgradient-type algorithms. In this paper, concrete nonsmooth convex optimization problems are discussed, and the two proposed algorithms are numerically compared with the existing  subgradient methods to evaluate their effectiveness.

This paper is organized as follows. Section \ref{sec:2} gives the mathematical preliminaries. Section \ref{sec:3} presents the incremental proximal point algorithm based on the Halpern fixed point algorithm and analyzes its convergence. Section \ref{sec:4} presents the incremental proximal point algorithm based on the Krasnosel'ski\u\i-Mann fixed point algorithm and analyzes its convergence. Section \ref{sec:5} describes concrete nonsmooth convex optimization problems and numerically compares the behaviors of the two proposed algorithms with those of the existing algorithms. Section \ref{sec:6} concludes the paper with a brief summary
and mentions future directions for improving the proposed algorithms.

\section{Mathematical Preliminaries}\label{sec:2}
Let $H$ be a real Hilbert space with inner product $\langle \cdot, \cdot \rangle$ and its induced norm $\|\cdot\|$, 
let $\mathbb{R}$ be the set of all real numbers, and let $\mathbb{N}$ be the set of all positive integers including zero.
Let $\mathrm{Fix}(T) := \{ x\in H \colon T(x) =x\}$ be the fixed point set of a mapping $T\colon H \to H$.
Let $\mathrm{dom}(f) := \{ x\in H \colon f(x) < \infty \}$ be the domain of a function $f \colon H \to (-\infty,\infty]$.
The identity mapping on $H$ is denoted by $\mathrm{Id}$.
Let $(x_n)_{n\in \mathbb{N}}$ be a sequence in $H$. 
A point $x\in H$ is said to be a {\em weak sequential cluster point} of $(x_n)_{n\in \mathbb{N}}$ \cite[subchapters 1.7 and 2.5]{b-c} if $(x_n)_{n\in \mathbb{N}}$ possesses a subsequence that weakly converges to $x\in H$.

\subsection{Nonexpansive mappings and proximity operators}
A mapping $T \colon H \to H$ is said to be {\em nonexpansive} \cite[Definition 4.1(ii)]{b-c} if $\|T(x) - T(y)\|\leq \|x-y\|$ $(x,y\in H)$.
$T$ is said to be {\em firmly nonexpansive} \cite[Definition 4.1(i)]{b-c} if 
$\|T(x) - T(y)\|^2 + \| (\mathrm{Id} - T)(x) - (\mathrm{Id} - T)(y) \|^2 \leq \|x-y\|^2$ $(x,y\in H)$.
The {\em metric projection} $P_C$ onto a nonempty, closed convex subset $C$ of $H$ 
is firmly nonexpansive with $\mathrm{Fix}(P_C) = C$ \cite[Proposition 4.8, (4.8)]{b-c}.

Let $f\colon H \to (-\infty,\infty]$ be proper, lower semicontinuous, and convex.
Then, the {\em proximity operator} of $f$ \cite[Definition 12.23]{b-c}, \cite{moreau1962}, denoted by $\mathrm{Prox}_f$, maps every $x\in H$ to the unique minimizer of 
$f + (1/2) \| x - \cdot \|^2$; i.e.,
\begin{align*}
\mathrm{Prox}_f (x)  =  \argmin_{y\in H} \left[ f(y) + \frac{1}{2} \left\|x-y\right\|^2  \right]
\text{ } \left(x\in H \right).
\end{align*}
The uniqueness and existence of $\mathrm{Prox}_f(x)$ are guaranteed for all $x\in H$ \cite[Definition 12.23]{b-c}, \cite{minty1965}.
The {\em subdifferential} of $f$ is the set-valued operator 
\begin{align*}
\partial f \colon H \to 2^H \colon x \mapsto \left\{ u\in H \colon f(y) \geq f(x) + \left\langle y-x,u \right\rangle \text{ } (y\in H)   \right\}.
\end{align*}

\begin{prop}{\em \cite[Propositions 12.26, 12.27, 12.28, and 16.14]{b-c}}\label{prop:1}
Let $f\colon H \to (-\infty,\infty]$ be proper, lower semicontinuous, and convex.
Then, the following hold:
\begin{enumerate}
\item[{\em(i)}]
Let $x,p\in H$. $p = \mathrm{Prox}_f (x)$ if and only if $x-p \in \partial f (p)$
(i.e., $\langle y-p,x-p  \rangle + f(p) \leq f(y)$ for all $y\in H$).
\item[{\em(ii)}]
$\mathrm{Prox}_f$ is firmly nonexpansive with $\mathrm{Fix}(\mathrm{Prox}_f) = \argmin_{x\in H}f(x)$.
\item[{\em (iii)}]
If $f$ is continuous at $x\in \mathrm{dom}(f)$,
$\partial f(x)$ is nonempty. Moreover, $\delta > 0$ exists such that $\partial f (B(x;\delta))$ is bounded, where $B(x;\delta)$ stands for a closed ball with center $x$
and radius $\delta$.
\end{enumerate}
\end{prop}

\subsection{Nonsmooth convex optimization problem with fixed point constraints}
Consider a networked system consisting of $I$ users, where user $i$ ($i\in \mathcal{I} := \{1,2,\ldots,I\}$) is assumed to have its own private mappings\footnote{The explicit forms of $T^{(i)}$ and $f^{(i)}$ are user $i$'s private information; i.e., other users cannot get the explicit forms of $T^{(i)}$ and $f^{(i)}$.} $T^{(i)}$ and $f^{(i)}$. The following problem is discussed.
\begin{prob}\label{problem:1}
Assume that 
\begin{enumerate}
\item[{\em (A1)}]
$T^{(i)} \colon H \to H$ ($i\in \mathcal{I}$) is firmly nonexpansive with 
$\bigcap_{i\in \mathcal{I}} \mathrm{Fix}(T^{(i)}) \neq \emptyset$;
\item[{\em (A2)}]
$f^{(i)} \colon H \to \mathbb{R}$ ($i\in \mathcal{I}$) is continuous and convex with $\mathrm{dom}(f^{(i)}) = H$ and
$\mathrm{Prox}_{f^{(i)}}$ can be efficiently computed.\footnote{Tables 10.1 and 10.2 in \cite{comb2011} present important examples of convex functions for which proximity operators can be easily computed within
a finite number of arithmetic operations.}
\end{enumerate}
Then,
\begin{align*}
\text{minimize } f(x) := \sum_{i\in \mathcal{I}} f^{(i)} (x) \text{ subject  to } x \in X := \bigcap_{i\in \mathcal{I}} \mathrm{Fix} \left(T^{(i)} \right).
\end{align*}
\end{prob}
The existence of a solution to Problem \ref{problem:1} is guaranteed when at least one of $\mathrm{Fix} (T^{(i)} )$ is bounded \cite[Theorem 25.C]{zeidler2b}.
Under the assumptions in the main theorems (Theorems \ref{theorem:2} and \ref{theorem:KM}), the existence of a solution to Problem \ref{problem:1} is guaranteed 
(see Lemmas \ref{lem:3-1}, \ref{lem:3}(iv), and \ref{lem:4}(iv)). 
If at least one of $f^{(i)}$ is strictly convex, the uniqueness of the solution to Problem \ref{problem:1} is also guaranteed \cite[Corollary 25.15]{zeidler2b}.

The following propositions will be used to prove the main theorems in this paper.

\begin{prop}{\em \cite[Lemma 1.2]{berinde}}\label{xu}
Assume that $(a_n)_{n\in \mathbb{N}} \subset [0,\infty)$ satisfies  
$a_{n+1} \leq (1- \alpha_n) a_n + \alpha_n \beta_n$ $(n\in \mathbb{N})$,
where $(\alpha_n)_{n\in \mathbb{N}} \subset (0,1]$ and $(\beta_n)_{n\in \mathbb{N}} \subset \mathbb{R}$ with
$\sum_{n = 1}^{\infty} \alpha_n = \infty$ and $\limsup_{n \to \infty} \beta_n \leq 0$.
Then, $\lim_{n \to \infty} a_n = 0$.
\end{prop}

\begin{prop}{\em \cite[Lemma 3.1]{opial}}\label{opial}
Suppose that $(x_n)_{n\in\mathbb{N}} \subset H$ weakly converges to $\hat{x} \in H$
and $\bar{x} \neq \hat{x}$.
Then,
$\liminf_{n\to\infty} \| x_n - \hat{x} \| < \liminf_{n\to\infty} \| x_n - \bar{x} \|$.
\end{prop}

\begin{prop}{\em \cite[Theorem 9.1]{b-c}}\label{bc}
When $f \colon H \to \mathbb{R}$ is convex, $f$ is weakly lower semicontinuous if and only if $f$ is lower semicontinuous. 
\end{prop}

\section{Halpern-type Incremental Proximal Point Algorithm}\label{sec:3}
This section presents the following algorithm for solving Problem \ref{problem:1}
using the Halpern algorithm \cite{halpern,wit} for finding a fixed point of a nonexpansive mapping.

\begin{algo}\label{algorithm:2}
\text{ }

\begin{enumerate}
\item[Step 0.] 
User $i$ ($i\in \mathcal{I}$) chooses $x^{(i)} \in H$ arbitrarily and sets 
$(\alpha_n)_{n\in \mathbb{N}} \subset (0,1]$ and $(\gamma_n)_{n\in\mathbb{N}} \subset (0,\infty)$.
User $I$ sets $x_0 \in H$ arbitrarily and transmits $x_0^{(0)} := x_0 \in H$ to user $1$.
\item[Step 1.]
User $i$ ($i\in \mathcal{I}$) computes $x_{n}^{(i)} \in H$ cyclically using
\begin{align*}
x_{n}^{(i)} := \alpha_n x^{(i)} + (1 - \alpha_n) T^{(i)} \left( \mathrm{Prox}_{\gamma_n f^{(i)}} \left( x_n^{(i-1)} \right) \right) \text{ } 
(i=1,2,\ldots,I).
\end{align*}
\item[Step 2.]
User $I$ defines $x_{n+1} \in H$ using
$x_{n+1} := x_n^{(I)}$
and transmits $x_{n+1}^{(0)} := x_{n+1}$ to user $1$.
The value of $n$ is then set to $n+1$, and the processing returns to Step 1.
\end{enumerate}
\end{algo}

The stopping criterions of Algorithm \ref{algorithm:2} are given by, for example, $\sum_{i\in \mathcal{I}} \| x_n - T^{(i)} (x_n) \| < \epsilon_1$ and $| f(x_{n-1}) - f(x_n) | < \epsilon_2$, where $\epsilon_i > 0$ $(i=1,2)$ is small enough. However, in general, such stopping criterions cannot be included in Algorithm \ref{algorithm:2} because none of the users can use all $x_n$, all $T^{(i)}$, and all $f^{(i)}$. If there exists an operator who manages the networked system and communicates with all users, the operator can verify whether the stopping criterions of Algorithm \ref{algorithm:2} are satisfied. The numerical section provides the number of iterations and elapsed time such that Algorithm \ref{algorithm:2} satisfies $| f(x_{n-1}) - f(x_n) | < 10^{-3}$ (see section \ref{sec:5} for details).

All users participating in the network are assumed to have the following information before the algorithm is executed.

\begin{assum}\label{assumption:2-1}
User $i$ ($i\in \mathcal{I}$) uses $(\alpha_n)_{n\in \mathbb{N}} \subset (0,1]$
and $(\gamma_n)_{n\in\mathbb{N}} \subset (0,\infty)$, 
which converge to $0$ and satisfy the following conditions:\footnote{Examples of $(\gamma_n)_{n\in \mathbb{N}}$ 
and $(\alpha_n)_{n\in\mathbb{N}}$ are 
$\gamma_n := 1/(n+1)^a$ and $\alpha_n := 1/(n+1)^b$ $(a\in (0,1/2), b \in (a, 1-a), a+b <1)$.}
\begin{align*}
&\text{{\em (C1)}}  \sum_{n=0}^{\infty} \alpha_n = \infty, \text{ } 
\text{{\em (C2)}} \lim_{n \to \infty} \frac{1}{\alpha_{n+1}} \left|\frac{1}{\gamma_{n+1}} - \frac{1}{\gamma_n} \right| = 0, \text{ } 
\text{{\em (C3)}} \lim_{n \to \infty} \frac{1}{\gamma_{n+1}} \left|1 - \frac{\alpha_n}{\alpha_{n+1}} \right| = 0,\\
&\text{{\em (C4) }} \lim_{n \to \infty} \frac{1}{\alpha_{n+1}}  \frac{\left| \gamma_{n+1} - \gamma_n\right|}{\gamma_{n+1}^2}  = 0, \text{ }
\text{{\em (C5)}} \lim_{n \to \infty} \frac{\alpha_n}{\gamma_n} = 0.
\end{align*}
\end{assum}

Step 1 in Algorithm \ref{algorithm:2} is a search for the fixed point of $T^{(i)}$, which is based on the Halpern algorithm \cite{halpern,wit} defined by $x_0 \in H$ and $x_{n+1} = \alpha_n x_0 + (1 - \alpha_n) T^{(i)} (x_n)$ $(n\in \mathbb{N})$. The algorithm with $\lim_{n\to\infty}\alpha_n = 0$ and (C1) strongly converges to the minimizer of $\| \cdot - x_0 \|^2$ over $\mathrm{Fix}(T^{(i)})$ \cite{halpern,wit}. Moreover, since $x_n^{(i)}$ in Step 1 uses the proximity operator $\mathrm{Prox}_{\gamma_n f^{(i)}}$, it can be seen intuitively that, for all $i\in \mathcal{I}$, $(x_n^{(i)})_{n\in \mathbb{N}}$ in Step 1 converges to not only a fixed point of $T^{(i)}$ but also a minimizer of $f^{(i)}$. Furthermore, Steps 1 and 2 in Algorithm \ref{algorithm:2} lead to the finding that $x_{n+1} = x_n^{(I)}(x_n^{(I-1)}) = x_n^{(I)}(x_n^{(I-1)}(x_n^{(I-2)})) = \cdots = x_n^{(I)}(x_n^{(I-1)},x_n^{(I-2)},\ldots, x_n^{(1)})$; i.e., $x_{n+1}$ has all the information of $x_n^{(i)}$ $(i\in \mathcal{I})$ needed to optimize $f^{(i)}$ over $\mathrm{Fix}(T^{(i)})$. Hence, it can be seen that $(x_n)_{n\in \mathbb{N}}$ approximates a minimizer of $\sum_{i\in \mathcal{I}} f^{(i)}$ over $\bigcap_{i\in \mathcal{I}} \mathrm{Fix}(T^{(i)})$. See subsection \ref{subsec:3.1} for the proof for the convergence property of $(x_n)_{n\in \mathbb{N}}$ in Algorithm \ref{algorithm:2}.
%

This convergence result depends on the following assumption.

\begin{assum}\label{assumption:3-1}
The sequence $(y_n^{(i)} := \mathrm{Prox}_{\gamma_n f^{(i)}} (x_n^{(i-1)}) )_{n\in \mathbb{N}}$ ($i\in \mathcal{I}$) 
generated by Algorithm \ref{algorithm:2} is bounded.
\end{assum}

Assume that, for all $i\in\mathcal{I}$, $\argmin_{x\in H} f^{(i)}(x) (= \mathrm{Fix}(\mathrm{Prox}_{f^{(i)}}) ) \neq \emptyset$ and $\mathrm{Fix}(T^{(i)})$ is bounded.
Then, user $i$ can choose in advance a bounded, closed convex set $X^{(i)}$ (e.g., $X^{(i)}$ is a closed ball with a large enough radius) 
satisfying $X^{(i)} \supset \mathrm{Fix}(T^{(i)})$.
Accordingly, user $i$ can compute, for example, 
\begin{align}\label{y_n_2}
x_{n}^{(i)} := P_{X^{(i)}} \left[ \alpha_n x^{(i)} + (1 - \alpha_n) T^{(i)} \left(y_n^{(i)} \right) \right]
\end{align}
instead of $x_n^{(i)}$ in Algorithm \ref{algorithm:2}.
Since $X^{(i)}$ $(i\in\mathcal{I})$ is bounded, 
$(x_n^{(i)})_{n\in\mathbb{N}}$ $(i\in \mathcal{I})$ is bounded.
Moreover, since Proposition \ref{prop:1}(ii) ensures that 
$\| y_n^{(i)} -x \| \leq \| x_n^{(i-1)} -x \|$ ($i\in \mathcal{I}, x\in \mathrm{Fix}(\mathrm{Prox}_{f^{(i)}})$),
the boundedness of $(x_n^{(i)})_{n\in\mathbb{N}}$ ($i\in\mathcal{I}$) guarantees that $(y_n^{(i)})_{n\in\mathbb{N}}$ ($i\in \mathcal{I}$) is bounded.
Hence, it can be assumed that $(x_n^{(i)})_{n\in \mathbb{N}}$ $(i\in \mathcal{I})$ in Algorithm \ref{algorithm:2} is as in \eqref{y_n_2}
in place of Assumption \ref{assumption:3-1}.

Next, a convergence analysis of Algorithm \ref{algorithm:2} is presented.

\begin{thm}\label{theorem:2}
Under Assumptions (A1), (A2), \ref{assumption:2-1}, and \ref{assumption:3-1}, any weak sequential cluster point of $(x_n^{(i)})_{n\in\mathbb{N}}$ ($i\in \mathcal{I}$) generated by Algorithm \ref{algorithm:2} belongs to the solution set of Problem \ref{problem:1}.
\end{thm}

An application example of Algorithm \ref{algorithm:2} is as follows. 
Let $X^{(i)} \subset H$ $(i\in \mathcal{I})$ be bounded, closed, and convex 
(see \eqref{y_n_2}),
let $C_k^{(i)} \subset H$ $(i\in \mathcal{I}, k\in \mathcal{K}^{(i)} := \{1,2,\ldots,K^{(i)}\})$ be a closed convex set onto which the projection can be easily calculated,
and let $(w_k^{(i)})_{k\in \mathcal{K}^{(i)}} \subset (0,1)$ $(i\in \mathcal{I})$ satisfy 
$\sum_{k\in \mathcal{K}^{(i)}} w_k^{(i)} = 1$.
Here let us define 
\begin{align}
&g^{(i)}(x) := \frac{1}{2} \sum_{k\in \mathcal{K}^{(i)}} w_k^{(i)}
\left( \min_{y\in C_k^{(i)}} \| x - y \| \right)^2 \text{ } (x\in H),\label{g}\\
&\bar{T}^{(i)} := P_{X^{(i)}} \left(\sum_{k\in \mathcal{K}^{(i)}} w_k^{(i)} P_{C_k^{(i)}} \right),
\text{ and } 
T^{(i)} := \frac{1}{2} \left(\mathrm{Id} + \bar{T}^{(i)} \right)
\text{ } (i\in \mathcal{I}).\label{T}
\end{align}
The function $g^{(i)}$ $(i\in \mathcal{I})$ defined by \eqref{g}
stands for the mean square value of the distances from $x\in H$ to $C_k^{(i)}$s.
Accordingly, we can express a subset of $X^{(i)}$ with the elements closest to $C_k^{(i)}$s in terms of the mean square norm by 
\begin{align}\label{gcfs}
C_{g^{(i)}} := \left\{ x \in X^{(i)} \colon g^{(i)}(x) = \min_{y\in X^{(i)}}
g^{(i)}(y) \right\} \text{ } (i\in \mathcal{I}). 
\end{align}
The $C_{g^{(i)}}$ is referred to as the {\em generalized convex feasible set} \cite[section I, Framework 2]{com1999}, \cite[Definition 4.1]{yamada}.
The condition $C_{g^{(i)}} \neq \emptyset$ $(i\in \mathcal{I})$ holds
from the boundedness of $X^{(i)}$ \cite[Remark 4.3(a)]{yamada}.
Even if $X^{(i)} \cap \bigcap_{k\in \mathcal{K}^{(i)}} C_k^{(i)} = \emptyset$, $C_{g^{(i)}}$ is well-defined.
In particular, $C_{g^{(i)}} = X^{(i)} \cap \bigcap_{k\in \mathcal{K}^{(i)}} C_k^{(i)}$ holds when $X^{(i)} \cap \bigcap_{k\in \mathcal{K}^{(i)}} C_k^{(i)} \neq \emptyset$.
Furthermore, $\bar{T}^{(i)}$ $(i\in \mathcal{I})$ is nonexpansive with $\mathrm{Fix}(\bar{T}^{(i)}) = C_{g^{(i)}}$ \cite[Proposition 4.2]{yamada}; i.e., 
$T^{(i)}$ $(i\in \mathcal{I})$ defined by \eqref{T} is firmly nonexpansive and 
\begin{align*}
\mathrm{Fix}\left(T^{(i)}\right) = \mathrm{Fix}\left(\bar{T}^{(i)}\right)
= C_{g^{(i)}} \text{ } (i\in \mathcal{I}). 
\end{align*}

Therefore, Theorem \ref{theorem:2} leads to the following.
\begin{cor}\label{cor:1}
Let $T^{(i)}$ ($i\in \mathcal{I}$) be a mapping defined by \eqref{T},
let $f^{(i)}$ ($i\in \mathcal{I}$) satisfy (A2),
and let $(x_n^{(i)})_{n\in \mathbb{N}}$ ($i\in \mathcal{I}$) 
be the sequence generated by \eqref{y_n_2}, where $(\alpha_n)_{n\in \mathbb{N}}$ and 
$(\gamma_n)_{n\in \mathbb{N}}$ satisfy Assumption \ref{assumption:2-1}.
Then, any weak sequential cluster point of $(x_n^{(i)})_{n\in \mathbb{N}}$
($i\in \mathcal{I}$) belongs to the solution set of Problem \ref{problem:1} with 
$X = \bigcap_{i\in \mathcal{I}} C_{g^{(i)}}$.
\end{cor}
Section \ref{sec:5} applies the proposed algorithms to the problem 
of minimizing $f$ over $\bigcap_{i\in \mathcal{I}} C_{g^{(i)}}$
and compares the behaviors of the proposed algorithms with the existing ones.

\subsection{Proof of Theorem \ref{theorem:2}}\label{subsec:3.1}
First, the following lemma is proven.

\begin{lem}\label{lem:1-1}
Suppose that Assumptions (A1), (A2), and \ref{assumption:3-1} hold and 
$(x_n^{(i)})_{n\in\mathbb{N}}$ and $(y_n^{(i)})_{n\in\mathbb{N}}$ ($i\in\mathcal{I}$) are the sequences generated by Algorithm \ref{algorithm:2}.
Then,
$(T^{(i)} (y_n^{(i)}))_{n\in \mathbb{N}}$ and
$({x}_n^{(i)})_{n\in \mathbb{N}}$ ($i\in\mathcal{I}$)
are bounded. 
\end{lem}

{\em Proof.}
Assumption (A1) guarantees that $\| T^{(i)} (y_n^{(i)}) - x \| \leq \|y_n^{(i)} -x\|$ $(i\in \mathcal{I}, n\in\mathbb{N}, x\in X)$,
which, together with Assumption \ref{assumption:3-1}, implies that $(T^{(i)} (y_n^{(i)}))_{n\in \mathbb{N}}$ $(i\in\mathcal{I})$ is bounded.
The definition of $x_n^{(i)}$ $(i\in\mathcal{I},n\in\mathbb{N})$ and the boundedness of $(T^{(i)} (y_n^{(i)}))_{n\in \mathbb{N}}$
lead to the boundedness of $(x_n^{(i)})_{n\in\mathbb{N}}$ $(i\in \mathcal{I})$.
\qed

Next, the following lemma is considered.

\begin{lem}\label{lem:2-1}
Suppose that Assumptions (A1), (A2), \ref{assumption:2-1}, and \ref{assumption:3-1} are satisfied.
Then, the following hold:
\begin{enumerate}
\item[{\em(i)}]
$\lim_{n\to\infty} \|{x}_{n+1} - {x}_n\|/\gamma_n=0$;
\item[{\em(ii)}]
$\lim_{n\to\infty} \|{y}_n^{(i)} - T^{(i)} ({y}_n^{(i)}) \| =0$ and 
$\lim_{n\to\infty} \|{x}_n^{(i-1)} - {y}_n^{(i)} \| =0$  ($i\in \mathcal{I}$);
\item[{\em(iii)}]
$\lim_{n\to\infty} \|{x}_n - {y}_n^{(i)} \|=0$ and
$\lim_{n\to\infty} \|{x}_n - T^{(i)} ({x}_n) \|=0$ ($i\in \mathcal{I}$).
\end{enumerate}
\end{lem}

{\em Proof.}
(i)
The definition of $x_n^{(i)}$ $(i\in \mathcal{I}, n\in\mathbb{N})$ and (A1) imply that, for all $i\in \mathcal{I}$ and for all $n\geq 1$,
\begin{align*}
&\left\|{x}_{n}^{(i)} - {x}_{n-1}^{(i)} \right\|\\
=& \left\| (1-\alpha_n) \left(T^{(i)} \left({y}_n^{(i)}\right) - T^{(i)} \left({y}_{n-1}^{(i)} \right) \right) 
       + (\alpha_n - \alpha_{n-1}) \left({x}^{(i)} - T^{(i)} \left({y}_{n-1}^{(i)} \right) \right) \right\|\\
\leq& (1-\alpha_n) \left\|T^{(i)} \left( {y}_n^{(i)} \right) - T^{(i)} \left({y}_{n-1}^{(i)} \right)\right\|
       + |\alpha_n - \alpha_{n-1}| \left\| {x}^{(i)} - T^{(i)} \left({y}_{n-1}^{(i)} \right) \right\|\\
\leq& (1-\alpha_n) \left\| {y}_n^{(i)}  - {y}_{n-1}^{(i)} \right\|
       + M_1 |\alpha_n - \alpha_{n-1}|,    
\end{align*}
where $M_1 := \max_{i\in \mathcal{I}} (\sup\{ \| {x}^{(i)} - T^{(i)}({y}_{n}^{(i)}) \| \colon n\in \mathbb{N}  \})$ 
and $M_1 < \infty$ holds from Lemma \ref{lem:1-1}.
Given definition $\bar{y}_n^{(i)} := \mathrm{Prox}_{\gamma_{n+1} f^{(i)}} (x_n^{(i-1)})$ $(i\in \mathcal{I}, n\in\mathbb{N})$,
Proposition \ref{prop:1}(ii) ensures that, for all $i\in \mathcal{I}$ and for all $n\geq 1$, 
\begin{align*}
\left\|{y}_{n}^{(i)} - {y}_{n-1}^{(i)} \right\|
&\leq \left\| \mathrm{Prox}_{\gamma_n f^{(i)}} \left( x_n^{(i-1)}  \right) - \mathrm{Prox}_{\gamma_{n} f^{(i)}} \left( x_{n-1}^{(i-1)}  \right) \right\|
       + \left\| \bar{y}_{n-1}^{(i)}  - y_{n-1}^{(i)}  \right\|\\
&\leq \left\|  x_n^{(i-1)}  - x_{n-1}^{(i-1)}   \right\|
       + \left\| \bar{y}_{n-1}^{(i)}  - y_{n-1}^{(i)}  \right\|.       
\end{align*}
Proposition \ref{prop:1}(i) means that $y_{n-1}^{(i)} := \mathrm{Prox}_{\gamma_{n-1} f^{(i)}}(x_{n-1}^{(i-1)})$ 
and $\bar{y}_{n-1}^{(i)} := \mathrm{Prox}_{\gamma_{n} f^{(i)}} (x_{n-1}^{(i-1)})$ satisfy
$(x_{n-1}^{(i-1)} - y_{n-1}^{(i)})/\gamma_{n-1} \in \partial f^{(i)} (y_{n-1}^{(i)})$
and 
$(x_{n-1}^{(i-1)} - \bar{y}_{n-1}^{(i)})/\gamma_{n} \in \partial f^{(i)} (\bar{y}_{n-1}^{(i)})$.
Accordingly, the monotonicity of $\partial f^{(i)}$ guarantees that, for all $i\in\mathcal{I}$ and for all $n\geq 1$, 
\begin{align*}
\left\langle y_{n-1}^{(i)} - \bar{y}_{n-1}^{(i)}, 
\frac{x_{n-1}^{(i-1)} - y_{n-1}^{(i)}}{\gamma_{n-1}} - \frac{x_{n-1}^{(i-1)} - \bar{y}_{n-1}^{(i)}}{\gamma_n}   \right\rangle \geq 0.
\end{align*}
Hence, 
\begin{align*}
&\frac{1}{\gamma_{n-1} \gamma_n} \bigg\{
\left\langle y_{n-1}^{(i)} - \bar{y}_{n-1}^{(i)}, 
\left(\gamma_n - \gamma_{n-1} \right) x_{n-1}^{(i-1)} \right\rangle
+
\left\langle y_{n-1}^{(i)} - \bar{y}_{n-1}^{(i)}, 
- \gamma_n \left(y_{n-1}^{(i)} - \bar{y}_{n-1}^{(i)} \right) \right\rangle\\
&\quad +
\left\langle y_{n-1}^{(i)} - \bar{y}_{n-1}^{(i)},
\left(\gamma_{n-1} - \gamma_n \right) \bar{y}_{n-1}^{(i)} \right\rangle
\bigg\} \geq 0, 
\end{align*}
which, together with the triangle inequality, means that
\begin{align*}
\left\| y_{n-1}^{(i)} - \bar{y}_{n-1}^{(i)} \right\|^2 
&\leq \frac{|\gamma_n - \gamma_{n-1}|}{\gamma_n} \left( \left\| x_{n-1}^{(i-1)} \right\| + \left\| \bar{y}_{n-1}^{(i)} \right\|  \right) 
\left\| y_{n-1}^{(i)} - \bar{y}_{n-1}^{(i)} \right\|\\
&\leq M_2 \frac{|\gamma_n - \gamma_{n-1}|}{\gamma_n} \left\| y_{n-1}^{(i)} - \bar{y}_{n-1}^{(i)} \right\|,
\end{align*}
where $M_2 := \max_{i\in \mathcal{I}} (\sup\{ \| {x}_{n}^{(i-1)} \| + \| \bar{y}_n^{(i)} \| \colon n\in \mathbb{N}  \})$ 
and $M_2 < \infty$ holds from Lemma \ref{lem:1-1}, Assumption \ref{assumption:3-1}, and Proposition \ref{prop:1}(ii).
Thus, for all $i\in\mathcal{I}$ and for all $n\geq 1$,
\begin{align*}
\left\| y_{n-1}^{(i)} - \bar{y}_{n-1}^{(i)} \right\| \leq M_2 \frac{|\gamma_n - \gamma_{n-1}|}{\gamma_n}.
\end{align*}
Therefore, for all $i\in\mathcal{I}$ and for all $n\geq 1$,
\begin{align}\label{1}
\left\|{x}_{n}^{(i)} - {x}_{n-1}^{(i)} \right\| 
&\leq (1-\alpha_n) \left\|  x_n^{(i-1)}  - x_{n-1}^{(i-1)}   \right\| + M_2 \frac{|\gamma_n - \gamma_{n-1}|}{\gamma_n}
+ M_1 |\alpha_n - \alpha_{n-1}|,      
\end{align}
which implies that, for all $n\geq 1$,
\begin{align*}
\| x_{n+1} - x_n \| 
&\leq (1-\alpha_n) \left\| {x}_{n}  - {x}_{n-1} \right\| + I M_1 |\alpha_n - \alpha_{n-1}| + I M_2 \frac{|\gamma_n - \gamma_{n-1}|}{\gamma_n}.
\end{align*}
Hence, for all $n \geq 1$, 
\begin{align*}
\frac{\|{x}_{n+1} - {x}_n\|}{\gamma_n}
&\leq  (1- \alpha_n) \frac{\| {x}_{n} - {x}_{n-1}\|}{\gamma_{n-1}}  
    + (1- \alpha_n) \left\{ \frac{\| {x}_{n} - {x}_{n-1}\|}{\gamma_{n}} - \frac{\| {x}_{n} - {x}_{n-1}\|}{\gamma_{n -1}}  \right\}\\
&\quad   + I {M}_1 \frac{|\alpha_n -\alpha_{n-1}|}{\gamma_n} + IM_2 \frac{|\gamma_n - \gamma_{n-1}|}{\gamma_n^2}\\
&\leq  (1- \alpha_n) \frac{\| {x}_{n} - {x}_{n-1}\|}{\gamma_{n-1}} 
+ I {M}_1 \frac{|\alpha_n -\alpha_{n-1}|}{\gamma_n} + IM_2 \frac{|\gamma_n - \gamma_{n-1}|}{\gamma_n^2}\\
&\quad + M_3 \left| \frac{1}{\gamma_n} - \frac{1}{\gamma_{n-1}}  \right|,
\end{align*}
where $M_3 := \sup \{\|{x}_{n+1} - {x}_n \| \colon n\in \mathbb{N} \} < \infty$.
This leads to the finding that 
\begin{align*}
\frac{\|{x}_{n+1} - {x}_n\|}{\gamma_n}
\leq (1- \alpha_n) \frac{\| {x}_{n} - {x}_{n-1}\|}{\gamma_{n-1}} 
+ \alpha_n X_n \text{ } (n\geq 1),
\end{align*}
where 
\begin{align*}
X_n := 
 I {M}_1 \frac{1}{\alpha_n} \frac{|\alpha_n -\alpha_{n-1}|}{\gamma_n} + IM_2 \frac{1}{\alpha_n}\frac{|\gamma_n - \gamma_{n-1}|}{\gamma_n^2}
  + M_3 \frac{1}{\alpha_n}\left| \frac{1}{\gamma_n} - \frac{1}{\gamma_{n-1}}  \right| \text{ } (n\geq 1).
\end{align*}
Proposition \ref{xu} and (C1), (C2), (C3), and (C4) ensure that
\begin{align}\label{xn-1}
\lim_{n \to \infty} \frac{\|{x}_{n+1} - {x}_n\|}{\gamma_n} =0.
\end{align}
Equation \eqref{xn-1} and $\lim_{n\to \infty} \gamma_n = 0$ imply that
$\lim_{n \to \infty} \|{x}_{n+1} - {x}_n\| =0$.

(ii)
The convexity of $\| \cdot \|^2$ and (A1) guarantee that, for all $x\in X$, for all $n\in \mathbb{N}$, and for all $i\in \mathcal{I}$,
\begin{align*}
\left\| x_{n}^{(i)} - x \right\|^2 
&\leq \alpha_n \left\| x^{(i)} - x \right\|^2 + (1-\alpha_n) \left\| T^{(i)} \left( y_n^{(i)} \right) - T^{(i)}(x) \right\|^2\\
&\leq \alpha_n \left\| x^{(i)} - x \right\|^2 +  
 \left\| y_n^{(i)} - x \right\|^2 - (1-\alpha_n) \left\| y_n^{(i)} -  T^{(i)} \left( y_n^{(i)} \right) \right\|^2.
\end{align*}
Proposition \ref{prop:1}(i) and $y_n^{(i)} := \mathrm{Prox}_{\gamma_n f^{(i)}}(x_n^{(i-1)})$ $(i\in\mathcal{I},n\in\mathbb{N})$ mean that, for all $x\in X$, for all $n\in \mathbb{N}$, and for all $i\in \mathcal{I}$, 
\begin{align*}
\left\langle x - y_n^{(i)},x_n^{(i-1)} - y_n^{(i)} \right\rangle \leq \gamma_n \left( f^{(i)} (x) - f^{(i)} \left(y_n^{(i)} \right)  \right).
\end{align*}
Moreover, from $\langle x,y \rangle = (1/2) (\| x \|^2 + \| y \|^2 - \|x-y\|^2)$ $(x,y\in H)$,
\begin{align*}
\left\langle x - y_n^{(i)},x_n^{(i-1)} - y_n^{(i)} \right\rangle 
= \frac{1}{2} \left( \left\|x - y_n^{(i)} \right\|^2 
  + \left\|x_n^{(i-1)} - y_n^{(i)} \right\|^2 - \left\|x- x_n^{(i-1)} \right\|^2 \right) 
\end{align*}
for all $x\in X$, for all $i\in\mathcal{I}$, and for all $n\in\mathbb{N}$.
Hence, for all $x\in X$, for all $n\in \mathbb{N}$, and for all $i\in \mathcal{I}$,
\begin{align}\label{y}
\left\|y_n^{(i)} -x \right\|^2 \leq \left\|x_n^{(i-1)} - x \right\|^2 - \left\|x_n^{(i-1)} - y_n^{(i)}\right\|^2
+ 2 \gamma_n \left( f^{(i)} (x) - f^{(i)} \left(y_n^{(i)} \right) \right).
\end{align}
Accordingly, setting $M_4 := \max_{i\in\mathcal{I}} \| x^{(i)} - x \|^2$ $(x\in X)$ leads to
\begin{align}\label{con}
\left\| x_{n}^{(i)} - x \right\|^2 
&\leq M_4 \alpha_n - (1-\alpha_n) \left\| y_n^{(i)} -  T^{(i)} \left( y_n^{(i)} \right) \right\|^2\\
&\quad + \left\|x_n^{(i-1)} - x \right\|^2 - \left\|x_n^{(i-1)} - y_n^{(i)}\right\|^2
+ 2 \gamma_n \left( f^{(i)} (x) - f^{(i)} \left(y_n^{(i)} \right) \right).\nonumber
\end{align}
Since Proposition \ref{prop:1}(iii) and (A2) ensure the existence of $z^{(i)} \in \partial f^{(i)}(x)$ and the boundedness of $\partial f^{(i)}(x)$ 
$(x \in X, i\in\mathcal{I})$,
the definition of $\partial f^{(i)}$ and Assumption \ref{assumption:3-1} imply that there exists $M_5 < \infty$ such that
$2(f^{(i)} (x) - f^{(i)} (y_n^{(i)})) \leq 2 \langle x - y_n^{(i)}, z^{(i)} \rangle 
\leq 2 \| x - y_n^{(i)} \| \|z^{(i)}\| \leq M_5$ $(x\in X, i\in \mathcal{I},n\in\mathbb{N})$.
Hence, for all $x \in X$ and for all $n\in \mathbb{N}$,
\begin{align*}
\| x_{n+1} - x \|^2 
&\leq I M_4 \alpha_n  - (1-\alpha_n) \sum_{i\in\mathcal{I}} \left\| y_n^{(i)} -  T^{(i)} \left( y_n^{(i)} \right) \right\|^2
+ \left\| x_n - x \right\|^2  + I M_5 \gamma_n\\ 
&\quad - \sum_{i\in\mathcal{I}} \left\| x_n^{(i-1)} - y_n^{(i)}\right\|^2.
\end{align*}
Since Lemma \ref{lem:1-1} means the existence of $M_6 < \infty$ such that,
for all $x \in X$ and for all $n\in\mathbb{N}$, 
\begin{align*}
\left\| x_n - x \right\|^2 - \left\| x_{n+1} - x \right\|^2 
&= \left( \left\| x_n - x \right\| - \left\| x_{n+1} - x \right\|\right)
   \left( \left\| x_n - x \right\| + \left\| x_{n+1} - x \right\|\right)\\ 
&\leq M_6 \left\| x_{n+1} - x_n \right\|,
\end{align*}  
we have that, for all $x\in X$ and for all $n\in \mathbb{N}$,
\begin{align*}
(1-\alpha_n) \sum_{i\in \mathcal{I}} \left\| y_n^{(i)} - T^{(i)} \left( y_n^{(i)} \right) \right\|^2 
&\leq I M_4 \alpha_n + I M_5 \gamma_n + M_6 \| x_{n+1} - x_n \|,\\
\sum_{i\in\mathcal{I}} \left\| x_n^{(i-1)} - y_n^{(i)}\right\|^2
&\leq I M_4 \alpha_n + I M_5 \gamma_n + M_6 \| x_{n+1} - x_n \|,
\end{align*}
which, together with $\lim_{n \to \infty} \|{x}_{n+1} - {x}_n\| =0$ and $\lim_{n\to\infty} \alpha_n = \lim_{n\to \infty} \gamma_n = 0$, implies that
\begin{align}\label{key0}
\lim_{n\to\infty} \left\| y_n^{(i)} - T^{(i)} \left( y_n^{(i)} \right) \right\| = 0 \text{ and } 
\lim_{n\to\infty} \left\| x_n^{(i-1)} - y_n^{(i)}\right\| = 0 \text{ } (i\in \mathcal{I}).
\end{align}

(iii)
From $\| x_n^{(i)} - T^{(i)} (y_n^{(i)}) \| = \alpha_n \| x^{(i)} - T^{(i)} (y_n^{(i)}) \|$ 
$(i\in \mathcal{I}, n\in \mathbb{N})$ and $\lim_{n\to\infty} \alpha_n = 0$, 
$\lim_{n\to\infty} \| x_n^{(i)} - T^{(i)}(y_n^{(i)}) \| = 0$ $(i\in \mathcal{I})$.
Since, for all $i\in \mathcal{I}$ and for all $n\in \mathbb{N}$,
\begin{align*}
\left\| x_n - x_n^{(i-1)} \right\| 
&\leq 
\sum_{j=1}^{i-1}
\left( \left\| x_n^{(j-1)} - y_n^{(j)} \right\| + \left\| y_n^{(j)} - T^{(j)} (y_n^{(j)}) \right\| + \left\| T^{(j)} (y_n^{(j)}) - x_n^{(j)} \right\| \right),
\end{align*}
\eqref{key0} and $\lim_{n\to\infty} \| x_n^{(i)} - T^{(i)}(y_n^{(i)}) \| = 0$ $(i\in \mathcal{I})$ guarantee that
$\lim_{n\to\infty} \| x_n - x_n^{(i-1)} \| = 0$ $(i\in \mathcal{I})$.
From $\| y_n^{(i)} - x_n \| 
\leq \| y_n^{(i)}  - x_n^{(i-1)}\| + \| x_n^{(i-1)} -x_n \|$ $(i\in\mathcal{I},n\in\mathbb{N})$, 
\eqref{key0} implies that 
\begin{align}\label{key3}
\lim_{n\to\infty} \left\| x_n - y_n^{(i)} \right\|=0 \text{ } (i\in \mathcal{I}).
\end{align}
Moreover, since 
$\| x_n - T^{(i)} (x_n) \| \leq \| x_n -y_n^{(i)} \| 
+ \|y_n^{(i)} - T^{(i)} (y_n^{(i)})\| + \| T^{(i)} (y_n^{(i)}) - T^{(i)} (x_n) \|$ $(i\in\mathcal{I},n\in\mathbb{N})$, 
(A1), \eqref{key0}, and \eqref{key3} ensure that 
\begin{align}\label{Z}
\lim_{n\to\infty} \left\| x_n - T^{(i)} (x_n) \right\|= 0 \text{ } (i\in \mathcal{I}).
\end{align}
This proves Lemma \ref{lem:2-1}.
\qed

Lemmas \ref{lem:1-1} and \ref{lem:2-1} lead to the following lemma.

\begin{lem}\label{lem:3-1}
Suppose that the assumptions in Lemma \ref{lem:2-1} hold.
Then, the following hold:
\begin{enumerate}
\item[{\em (i)}]
$\limsup_{n\to\infty} f(x_n) \leq f(x)$ for all $x\in X$;
\item[{\em (ii)}]
There exists a weak sequential cluster point of $(x_n)_{n\in\mathbb{N}}$ that belongs to the solution set $X^\star$ of Problem \ref{problem:1};
\item[{\em (iii)}]
Any weak sequential cluster point of $(x_{n}^{(i)})_{n\in\mathbb{N}}$ ($i\in \mathcal{I}$) is in $X^\star$.
\end{enumerate}
\end{lem} 

{\em Proof.}
(i)
Inequality \eqref{con} guarantees that, for all $x\in X$, for all $n\in\mathbb{N}$, and for all $i\in\mathcal{I}$,
\begin{align*}
\left\|x_n^{(i)} - x \right\|^2
\leq \left\|x_n^{(i-1)} - x \right\|^2 + 2 \gamma_n \left( f^{(i)}(x) - f^{(i)} \left(y_n^{(i)} \right)  \right) + M_4 \alpha_n, 
\end{align*}
which, together with $x_{n+1} = x_n^{(I)} = x_{n+1}^{(0)}$ $(n\in\mathbb{N})$ and $f:= \sum_{i\in\mathcal{I}} f^{(i)}$, implies that 
\begin{align*}
\| x_{n+1} - x \|^2 
&\leq \| x_{n} - x \|^2 + 2 \gamma_n \sum_{i\in \mathcal{I}}\left( f^{(i)}(x) - f^{(i)} \left(y_n^{(i)} \right)  \right) + I M_4 \alpha_n\\
&= \| x_{n} - x \|^2 + I M_4 \alpha_n\\
&\quad + 2 \gamma_n \left( f (x) - f(x_n) + \sum_{i\in \mathcal{I}} \left[f^{(i)}(x_n) -  f^{(i)} \left(y_n^{(i)} \right) \right]  \right).
\end{align*}
Since Lemma \ref{lem:1-1} means that $M_6 < \infty$ exists such that $\| x_n - x \|^2 - \| x_{n+1} - x \|^2 \leq M_6 \| x_{n+1} - x_n \|$ 
$(x \in X, n\in\mathbb{N})$, 
for all $x\in X$ and for all $n\in\mathbb{N}$, 
\begin{align*}
2 \left(f(x_n) - f(x) \right) 
&\leq  \frac{M_6 \| x_n -  x_{n+1} \|}{\gamma_n} + I M_4 \frac{\alpha_n}{\gamma_n}
+ 2 \sum_{i\in \mathcal{I}} \left[ f^{(i)} (x_n) - f^{(i)} \left(y_n^{(i)} \right) \right].
\end{align*}
Moreover, the definition of $\partial f^{(i)}$ $(i\in \mathcal{I})$, (A2), Lemma \ref{lem:1-1}, and Proposition \ref{prop:1}(iii) lead to the existence of $M_7 < \infty$ such that, for all $i\in \mathcal{I}$ and for all $n\in \mathbb{N}$,
$f^{(i)} (x_n) - f^{(i)} (y_n^{(i)} ) \leq M_7 \| x_n - y_n^{(i)} \|$,
which, together with \eqref{key3}, implies that 
$\limsup_{n\to\infty} [ f^{(i)} (x_n) - f^{(i)} (y_n^{(i)} ) ] \leq 0$ $(i\in \mathcal{I})$.
Hence, \eqref{xn-1} and (C5) ensure that 
\begin{align*}
2 \limsup_{n\to\infty} \left(f(x_n) - f(x) \right) 
&\leq 2 \sum_{i\in \mathcal{I}} \limsup_{n\to\infty} \left[ f^{(i)} (x_n) - f^{(i)} \left(y_n^{(i)} \right) \right]
\leq 0.
\end{align*}
Therefore, $\limsup_{n\to\infty} f (x_n) \leq f (x)$ $(x\in X)$.

(ii)
Lemma \ref{lem:1-1} guarantees the existence of a weak sequential cluster point of $(x_n)_{n\in\mathbb{N}}$.
Let $x^* \in H$ be an arbitrary weak sequential cluster point of $(x_n)_{n\in\mathbb{N}}$.
Then, there exists 
$({x}_{n_k})_{k\in \mathbb{N}}$ $(\subset ({x}_n)_{n\in \mathbb{N}})$ such that $({x}_{n_k})_{k\in \mathbb{N}}$ weakly converges to $x^*$.
Here, $i\in \mathcal{I}$ is arbitrarily fixed, and $x^* \notin \mathrm{Fix}(T^{(i)})$ is assumed.
Then, Proposition \ref{opial}, Lemma \ref{lem:2-1}(iii), and (A1) produce a contradiction: 
\begin{align*}
\liminf_{k \to \infty} \left\|x_{n_k} - x^* \right\|  
&< \liminf_{k \to \infty} \left\| x_{n_k} -T^{(i)} (x^*) \right\|\\
&= \liminf_{k \to \infty} \left\|x_{n_k} -T^{(i)} \left(x_{n_k}\right) + T^{(i)} \left(x_{n_k}\right) - T^{(i)} (x^*) \right\|\\ 
&= \liminf_{k \to \infty} \left\| T^{(i)} \left(x_{n_k}\right) -T^{(i)} (x^*) \right\|\\
&\leq \liminf_{k \to \infty} \|x_{n_k} - {x}^*\|.
\end{align*}
Therefore, $x^* \in \mathrm{Fix}(T^{(i)})$ $(i\in \mathcal{I})$; i.e., $x^* \in X$.
Moreover, (A2), the weak convergence of $(x_{n_k})_{k\in \mathbb{N}}$ to $x^* \in X$, and 
Proposition \ref{bc} imply that $f(x^*) \leq \liminf_{k\to\infty} f(x_{n_k})$.
Accordingly, Lemma \ref{lem:3-1}(i) guarantees that, for all $x\in X$, 
\begin{align*}
f(x^*) \leq \liminf_{k\to\infty} f(x_{n_k}) \leq \limsup_{k\to\infty} f(x_{n_k}) 
\leq \limsup_{n\to\infty} f(x_{n})
\leq f(x); \text{ i.e., } x^* \in X^\star.
\end{align*}

(iii)
Lemma \ref{lem:3-1}(ii) means that any weak sequential cluster point of $(x_n)_{n\in\mathbb{N}}$ is in $X^\star$.
From Lemma \ref{lem:1-1}, $\lim_{n\to\infty} \| x_n - x_n^{(i-1)} \| = 0$ $(i\in\mathcal{I})$, and $x_{n+1} = x_n^{(I)}$ $(n\in\mathbb{N})$,
any weak sequential cluster point of $(x_n^{(i)})_{n\in\mathbb{N}}$ $(i\in\mathcal{I})$ is in $X^\star$.
This completes the proof.
\qed

\section{Krasnosel'ski\u\i-Mann-type Incremental Proximal Point Algorithm}\label{sec:4}
The following algorithm using the Krasnosel'ski\u\i-Mann algorithm \cite{kra,mann} 
is presented.

\begin{algo}\label{algorithm:KM}
\text{ }

\begin{enumerate}
\item[Step 0.] 
User $i$ ($i\in \mathcal{I}$) sets 
$(\alpha_n)_{n\in \mathbb{N}} \subset (0,1]$ and $(\gamma_n)_{n\in\mathbb{N}} \subset (0,\infty)$.
User $I$ sets $x_0 \in H$ arbitrarily and transmits $x_0^{(0)} := x_0 \in H$ to user $1$.
\item[Step 1.]
User $i$ ($i\in \mathcal{I}$) computes $x_{n}^{(i)} \in H$ cyclically using
\begin{align*}
x_{n}^{(i)} := \alpha_n x_n^{(i-1)} + (1 - \alpha_n) T^{(i)} \left( \mathrm{Prox}_{\gamma_n f^{(i)}} \left( x_n^{(i-1)} \right) \right) \text{ } 
(i=1,2,\ldots,I).
\end{align*}
\item[Step 2.]
User $I$ defines $x_{n+1} \in H$ using
$x_{n+1} := x_n^{(I)}$
and transmits $x_{n+1}^{(0)} := x_{n+1}$ to user $1$.
The value of $n$ is then set to $n+1$, and the processing returns to Step 1.
\end{enumerate}
\end{algo}

Two assumptions are made here.

\begin{assum}\label{assum:KM-1}
User $i$ ($i\in\mathcal{I}$) uses $(\alpha_n)_{n\in\mathbb{N}} \subset (0,1]$ and $(\gamma_n)_{n\in\mathbb{N}} \subset (0,\infty)$
satisfying the following conditions:\footnote{Examples of $(\gamma_n)_{n\in\mathbb{N}}$ and $(\alpha_n)_{n\in\mathbb{N}}$ are 
$\gamma_n := 1/(n+1)^a$ and $\alpha_n := t$ $(a \in (0,1], t \in (0,1))$.}
\begin{align*}
\text{{\em (C6) }} 0 < \liminf_{n\to \infty} \alpha_n \leq \limsup_{n\to \infty} \alpha_n < 1, 
\text{ {\em (C7) }} \lim_{n\to\infty} \gamma_n = 0, 
\text{ {\em (C8) }} \sum_{n=0}^\infty \gamma_n = \infty.
\end{align*}
\end{assum}
\begin{assum}\label{assum:KM-2}
The sequence $(y_n^{(i)} := \mathrm{Prox}_{\gamma_n f^{(i)}} (x_n^{(i-1)}) )_{n\in \mathbb{N}}$ ($i\in \mathcal{I}$) 
generated by Algorithm \ref{algorithm:KM} is bounded.
\end{assum}

Step 1 in Algorithm \ref{algorithm:KM} is a search for the fixed point of $T^{(i)}$, which is based on the Krasnosel'ski\u\i-Mann algorithm \cite{kra,mann} defined by $x_0\in H$ and $x_{n+1} = \alpha_n x_n + (1-\alpha_n) T^{(i)}(x_n)$ $(n\in \mathbb{N})$. It is guaranteed that the algorithm with (C6) weakly converges to a fixed point of $T^{(i)}$ \cite{kra,mann}. Accordingly, from the use of the proximity operator $\mathrm{Prox}_{\gamma_n f^{(i)}}$, it can be seen intuitively that $(x_n^{(i)})$ in Step 1 approximates a fixed point of $T^{(i)}$ as well as a minimizer of $f^{(i)}$. From the incremental steps in Steps 1 and 2 (see also the discussion of Algorithm \ref{algorithm:2}), it can be seen that Algorithm \ref{algorithm:KM} optimizes $\sum_{i\in \mathcal{I}} f^{(i)}$ over $\bigcap_{i\in \mathcal{I}} \mathrm{Fix}(T^{(i)})$. The mathematical proof for the convergence property of $(x_n)_{n\in \mathbb{N}}$ in Algorithm \ref{algorithm:KM} is given in subsection \ref{subsec:4.1}.
%

Next, a convergence analysis of Algorithm \ref{algorithm:KM} is presented.

\begin{thm}\label{theorem:KM}
Under Assumptions (A1), (A2), \ref{assum:KM-1}, and \ref{assum:KM-2}, 
there exists a weak sequential cluster point of $(x_n^{(i)})_{n\in\mathbb{N}}$ ($i\in \mathcal{I}$) generated by Algorithm \ref{algorithm:KM} 
which belongs to the solution set of Problem \ref{problem:1}.
\end{thm}

The discussion in section \ref{sec:3} leads to the following.
\begin{cor}\label{cor:2}
Let $T^{(i)}$ ($i\in \mathcal{I}$) be a mapping defined by \eqref{T},
let $f^{(i)}$ ($i\in \mathcal{I}$) satisfy (A2),
and let $(x_n^{(i)})_{n\in \mathbb{N}}$ ($i\in \mathcal{I}$) 
be the sequence generated by  
\eqref{y_n_2} when $x^{(i)}$ is replaced by $x_n^{(i-1)}$, where $(\alpha_n)_{n\in \mathbb{N}}$ and 
$(\gamma_n)_{n\in \mathbb{N}}$ satisfy Assumption \ref{assum:KM-1}.
Then, there exists a weak sequential cluster point of $(x_n^{(i)})_{n\in \mathbb{N}}$
($i\in \mathcal{I}$) which belongs to the solution set of Problem \ref{problem:1} with  
$X = \bigcap_{i\in \mathcal{I}} C_{g^{(i)}}$.
\end{cor}

\subsection{Proof of Theorem \ref{theorem:KM}}\label{subsec:4.1}
The proof starts with the following lemma.

\begin{lem}\label{lem:2}
The sequence $(x_n)_{n\in\mathbb{N}}$ generated by Algorithm \ref{algorithm:KM} satisfies that, for all $x\in X$ and for all $n\in\mathbb{N}$, 
\begin{align*}
\| x_{n+1} - x \|^2
&\leq \left\| x_n - x \right\|^2 - (1-\alpha_n) \sum_{i\in\mathcal{I}} \left\{ \left\|x_n^{(i-1)} - y_n^{(i)}\right\|^2
+  \left\| y_n^{(i)}  - T^{(i)} \left(y_n^{(i)} \right) \right\|^2 \right\}\nonumber\\
&\quad + 2 (1-\alpha_n)\gamma_n \sum_{i\in\mathcal{I}} \left[ f^{(i)} (x) - f^{(i)} \left(y_n^{(i)} \right) \right].
\end{align*}
\end{lem}

{\em Proof.}
The definition of $x_n^{(i)}$ $(i\in\mathcal{I},n\in\mathbb{N})$ and the convexity of $\| \cdot \|^2$ guarantee that,
for all $x\in X$, for all $n\in\mathbb{N}$, and for all $i\in\mathcal{I}$,
\begin{align*}
\left\| x_n^{(i)} - x \right\|^2
&\leq \alpha_n \left\| x_n^{(i-1)} - x \right\|^2 + (1-\alpha_n) \left\| T^{(i)} \left(y_n^{(i)} \right) - x \right\|^2,
\end{align*} 
which, together with (A1), implies that
\begin{align*}
\left\| x_n^{(i)} - x \right\|^2
&\leq \alpha_n \left\| x_n^{(i-1)} - x \right\|^2 + (1-\alpha_n) \left\| y_n^{(i)}  - x \right\|^2\\
&\quad - (1-\alpha_n) \left\| y_n^{(i)}  - T^{(i)} \left(y_n^{(i)} \right) \right\|^2.
\end{align*}
Moreover, \eqref{y} means that, for all $x\in X$, for all $n\in\mathbb{N}$, and for all $i\in\mathcal{I}$,
\begin{align*}
\left\| x_n^{(i)} - x \right\|^2
&\leq \alpha_n \left\| x_n^{(i-1)} - x \right\|^2 + (1-\alpha_n) \bigg\{
 \left\|x_n^{(i-1)} - x \right\|^2 - \left\|x_n^{(i-1)} - y_n^{(i)}\right\|^2\\
&\quad + 2 \gamma_n \left( f^{(i)} (x) - f^{(i)} \left(y_n^{(i)} \right) \right) \bigg\}
- (1-\alpha_n) \left\| y_n^{(i)}  - T^{(i)} \left(y_n^{(i)} \right) \right\|^2\\
&= \left\| x_n^{(i-1)} - x \right\|^2 - (1-\alpha_n)\left\|x_n^{(i-1)} - y_n^{(i)}\right\|^2\\
&\quad + 2 (1-\alpha_n)\gamma_n \left( f^{(i)} (x) - f^{(i)} \left(y_n^{(i)} \right) \right)
- (1-\alpha_n) \left\| y_n^{(i)}  - T^{(i)} \left(y_n^{(i)} \right) \right\|^2.
\end{align*}
Summing this inequality over all $i$ completes the proof of Lemma \ref{lem:2}.
\qed

The following lemma indicates that Theorem \ref{theorem:KM} holds when $(x_n)_{n\in\mathbb{N}}$ in Algorithm \ref{algorithm:KM} is Fej\'er monotone with respect to $X^\star$ \cite[chapter 5]{b-c}.

\begin{lem}\label{lem:3}
Suppose that Assumptions (A1), (A2), \ref{assum:KM-1}, and \ref{assum:KM-2} hold and there exists $n_0\in\mathbb{N}$ such that
$\| x_{n+1} - x^\star \| \leq \| x_n- x^\star \|$ for all $x^\star \in X^\star$ and for all $n\geq n_0$.
Then, the following hold:
\begin{enumerate}
\item[{\em (i)}]
$\lim_{n\to\infty} \| x_n^{(i-1)} - y_n^{(i)} \| = 0$ and $\lim_{n\to\infty} \| y_n^{(i)} - T^{(i)} (y_n^{(i)})  \|=0$
($i\in\mathcal{I}$);
\item[{\em (ii)}]
$\lim_{n\to\infty} \| x_n- y_n^{(i)} \| = 0$ and $\lim_{n\to\infty} \| x_n - T^{(i)} (x_n)  \|=0$
($i\in\mathcal{I}$);
\item[{\em (iii)}]
$\liminf_{n\to\infty} f(x_n) \leq f(x)$ ($x \in X$);
\item[{\em (iv)}]
There exists $(x_{n_l}^{(i)})_{l\in\mathbb{N}} \subset (x_n^{(i)})_{n\in\mathbb{N}}$ ($i\in\mathcal{I}$) which weakly converges to $x^* \in X^\star$.
\end{enumerate}
\end{lem}

{\em Proof.}
(i)
The definition of $\partial f^{(i)}$ ensures that, for all $x\in X$, for all $n\in\mathbb{N}$, and for all $i\in\mathcal{I}$,
$f^{(i)} (x) - f^{(i)} (y_n^{(i)}) \leq \langle x - y_n^{(i)}, z^{(i)} \rangle \leq N_1$, where
$z^{(i)} \in \partial f^{(i)}(x)$ $(i\in\mathcal{I})$,  
$N_1 := \max_{i\in\mathcal{I}}(\sup \{ \langle y_n^{(i)} - x, z^{(i)} \rangle \colon n\in \mathbb{N}\})$, and $N_1 < \infty$
is satisfied from Assumption \ref{assum:KM-2}.
Accordingly, Lemma \ref{lem:2} guarantees that, for all $x^\star \in X^\star$ and for all $n\in\mathbb{N}$,
\begin{align}\label{0}
\begin{split}
&(1-\alpha_n) \sum_{i\in\mathcal{I}}  \left\|x_n^{(i-1)} - y_n^{(i)}\right\|^2
\leq \left\| x_n - x^\star \right\|^2 - \| x_{n+1} - x^\star \|^2 + 2 I N_1 (1-\alpha_n)\gamma_n,\\
&(1-\alpha_n) \sum_{i\in\mathcal{I}}  \left\| y_n^{(i)}  - T^{(i)} \left(y_n^{(i)} \right) \right\|^2
\leq \left\| x_n - x^\star \right\|^2 - \| x_{n+1} - x^\star \|^2 + 2 I N_1 (1-\alpha_n)\gamma_n,
\end{split}
\end{align}
which, together with (C6), (C7), and the existence of $\lim_{n\to\infty} \|x_n -x^\star \|$  
(by $\| x_{n+1} - x^\star \| \leq \| x_n- x^\star \|$ ($x^\star \in X^\star, n\geq n_0$)), means that
$\lim_{n\to\infty} \|x_n^{(i-1)} - y_n^{(i)}\| = 0$ and 
$\lim_{n\to\infty} \| y_n^{(i)}  - T^{(i)} (y_n^{(i)} ) \| = 0$ $(i\in\mathcal{I})$.

(ii)
From 
\begin{align*}
\left\| x_n^{(i)} - x_n^{(i-1)} \right\|
&\leq  \left\| T^{(i)} \left(y_n^{(i)} \right) - x_n^{(i-1)}  \right\|
\leq \left\| T^{(i)} \left(y_n^{(i)} \right) -  y_n^{(i)} \right\| 
  + \left\| y_n^{(i)} - x_n^{(i-1)} \right\|,
\end{align*} 
Lemma \ref{lem:3}(i) leads to $\lim_{n\to\infty}\| x_n^{(i)} - x_n^{(i-1)} \| = 0$ $(i\in\mathcal{I})$.
Since 
\begin{align*}
\left\| T^{(i)} \left(y_n^{(i)} \right) - x_n^{(i)} \right\| 
\leq \left\|T^{(i)} \left(y_n^{(i)} \right) - y_n^{(i)} \right\|
+ \left\|y_n^{(i)} - x_n^{(i-1)} \right\| + \left\| x_n^{(i-1)} - x_n^{(i)} \right\|,
\end{align*}
Lemma \ref{lem:3}(i) implies that $\lim_{n\to\infty} \| T^{(i)}(y_n^{(i)}) - x_n^{(i)} \| = 0$ $(i\in\mathcal{I})$.
Thus, a discussion similar to the one for obtaining \eqref{key3} and \eqref{Z} leads to $\lim_{n\to\infty} \|x_n - x_n^{(i-1)}\| = 0$,
$\lim_{n\to\infty} \|x_n - y_n^{(i)}\| = 0$, 
and 
$\lim_{n\to\infty} \|x_n - T^{(i)} (x_n) \| = 0$ $(i\in\mathcal{I})$.

(iii)
From $f:= \sum_{i\in\mathcal{I}} f^{(i)}$, the definition of $\partial f^{(i)}$, (A2), and Proposition \ref{prop:1}(iii), 
there exists $N_2 < \infty$ such that, for all $x\in X$ and for all $n\in\mathbb{N}$,
\begin{align*}
\sum_{i\in\mathcal{I}} \left[ f^{(i)} (x) - f^{(i)} \left(y_n^{(i)} \right) \right]
&= f (x) - f(x_n) + \sum_{i\in\mathcal{I}} \left[ f^{(i)} (x_n) - f^{(i)} \left(y_n^{(i)} \right) \right]\\
&\leq f (x) - f(x_n) +  N_2  \sum_{i\in\mathcal{I}} \left\| x_n - y_n^{(i)} \right\|.
\end{align*}
Accordingly, Lemma \ref{lem:2} implies that, for all $x\in X$ and for all $n\in\mathbb{N}$,
\begin{align}\label{4.2}
2 (1-\alpha_n)\gamma_n \left(f (x_n) - f(x) - N_2  \sum_{i\in\mathcal{I}} \left\| x_n - y_n^{(i)} \right\| \right) \leq \left\| x_n - x \right\|^2 - \| x_{n+1} - x \|^2.
\end{align}
Summing up \eqref{4.2} from $n=0$ to infinity leads to
\begin{align*}
\sum_{n=0}^\infty \gamma_n (1-\alpha_n) \left(f (x_n) - f(x) - N_2  \sum_{i\in\mathcal{I}} \left\| x_n - y_n^{(i)} \right\| \right) \leq \| x_0 - x \| < \infty.
\end{align*}
It is next shown that $\liminf_{n\to\infty} (1-\alpha_n) (f (x_n) - f(x) - N_2  \sum_{i\in\mathcal{I}} \| x_n - y_n^{(i)}\| ) \leq 0$ $(x\in X)$.
If this assertion does not hold, there exist $m_0 \in \mathbb{N}$ and $\gamma > 0$ such that $(1-\alpha_n) (f (x_n) - f(x) - N_2  \sum_{i\in\mathcal{I}} \| x_n - y_n^{(i)}\| ) \geq \gamma$
for all $n \geq m_0$. 
Accordingly, (C8) ensures that, for all $x\in X$,
\begin{align*}
\infty = \gamma \sum_{n=m_0}^\infty \gamma_n \leq \sum_{n=m_0}^\infty \gamma_n (1-\alpha_n) \left(f (x_n) - f(x) - N_2  \sum_{i\in\mathcal{I}} \left\| x_n - y_n^{(i)} \right\| \right) < \infty,
\end{align*}
which is a contradiction. 
Hence, 
(C6) and Lemma \ref{lem:3}(ii) imply that there exists $\alpha \in (0,1)$ such that, for all $x\in X$,
\begin{align*}
(1-\alpha)  \liminf_{n\to\infty} \left( f (x_n) - f(x) \right) &\leq \liminf_{n\to\infty} (1-\alpha_n) \left( f (x_n) - f(x) \right)\\
&\leq N_2 \limsup_{n\to\infty}  (1-\alpha_n)  \sum_{i\in\mathcal{I}} \left\| x_n - y_n^{(i)} \right\|
= 0.
\end{align*}
Therefore, $\liminf_{n\to\infty} f(x_n) \leq f(x)$ $(x\in X)$.

(iv)
Lemma \ref{lem:3}(iii) ensures the existence of a subsequence $(x_{n_l})_{l\in\mathbb{N}}$ of $(x_n)_{n\in\mathbb{N}}$ such that,
for all $x\in X$,
\begin{align*}
\lim_{l\to\infty} f\left( x_{n_l}  \right) = \liminf_{n\to\infty} f(x_n) \leq f(x).
\end{align*}
The boundedness of $(x_{n_l})_{l\in\mathbb{N}}$ guarantees that there exists $(x_{n_{l_m}})_{m\in\mathbb{N}} \subset (x_{n_l})_{l\in\mathbb{N}}$
that weakly converges to $x^*$.
The same discussion as in the proof of Lemma \ref{lem:3-1}(ii) leads to $x^*\in X$. 
Since Proposition \ref{bc} implies that $f(x^* ) \leq \liminf_{m\to\infty} f(x_{n_{l_m}})$,
\begin{align*}
f(x^* ) \leq \liminf_{m\to\infty} f\left(x_{n_{l_m}}\right) = \lim_{l\to\infty} f \left( x_{n_l}\right) \leq f(x) \text{ } (x\in X), 
\text{ i.e., } x^* \in X^\star.
\end{align*}

Consider another subsequence $(x_{n_{l_k}})_{k\in\mathbb{N}}\subset (x_{n_l})_{l\in\mathbb{N}}$ that weakly converges to $x_*$.
From the above discussion, $x_* \in X^\star$.
Here, assume that $x^* \neq x_*$. 
Then, the existence of $\lim_{n\to\infty} \|x_n - x^\star \|$ $(x^\star \in X^\star)$ and Proposition \ref{opial} lead to a contradiction:
\begin{align*}
\lim_{n\to\infty} \left\| x_n - x^* \right\| 
&= \lim_{m\to\infty} \left\| x_{n_{l_m}} - x^* \right\| <  \lim_{m\to\infty} \left\| x_{n_{l_m}} - x_* \right\|\\
&= \lim_{n\to\infty} \left\| x_n - x_* \right\| = \lim_{k\to\infty} \left\| x_{n_{l_k}} - x_* \right\| < \lim_{k\to\infty} \left\| x_{n_{l_k}} - x^* \right\|\\
&= \lim_{n\to\infty} \left\| x_n - x^* \right\|.
\end{align*}
Therefore, any subsequence of $(x_{n_l})_{l\in\mathbb{N}}$ converges weakly to $x^* \in X^\star$.
This means that $(x_{n_l})_{l\in\mathbb{N}}$ weakly converges to $x^* \in X^\star$.
From $\lim_{n\to\infty} \| x_n - x_n^{(i-1)} \|=0$ $(i\in\mathcal{I})$, $(x_{n_l}^{(i)})_{n\in\mathbb{N}}$ $(i\in\mathcal{I})$ weakly converges to $x^* \in X^\star$.
This completes the proof.
\qed

Next it is proven that Theorem \ref{theorem:KM} holds when $(x_n)_{n\in\mathbb{N}}$ in Algorithm \ref{algorithm:KM} is not Fej\'er monotone with respect to $X^\star$.

\begin{lem}\label{lem:4}
Suppose that Assumptions (A1), (A2), \ref{assum:KM-1}, and \ref{assum:KM-2} hold and there exist
$x_0^\star \in X^\star$ and $(x_{n_j})_{j\in\mathbb{N}} \subset (x_n)_{n\in \mathbb{N}}$ such that
$\| x_{n_j} - x_0^\star \| < \| x_{n_j + 1}- x_0^\star \|$ for all $j\in \mathbb{N}$.
Then, the following hold:
\begin{enumerate}
\item[{\em (i)}]
$\lim_{j\to\infty} \| x_{n_j}^{(i-1)} - y_{n_j}^{(i)} \| = 0$ and $\lim_{j\to\infty} \| y_{n_j}^{(i)} - T^{(i)} (y_{n_j}^{(i)})  \|=0$
($i\in\mathcal{I}$);
\item[{\em (ii)}]
$\lim_{j\to\infty} \| x_{n_j}- y_{n_j}^{(i)} \| = 0$ and $\lim_{j\to\infty} \| x_{n_j} - T^{(i)} (x_{n_j})  \|=0$
($i\in\mathcal{I}$);
\item[{\em (iii)}]
$\limsup_{j\to\infty} f(x_{n_j}) \leq f(x_0^\star)$;
\item[{\em (iv)}]
There exists a weak sequential cluster point of $(x_n)_{n\in\mathbb{N}}$ which is in $X^\star$.
\end{enumerate}
\end{lem}

{\em Proof.}
(i)
A discussion similar to the one for obtaining \eqref{0} and $\| x_{n_j} - x_0^\star \| < \| x_{n_j + 1}- x_0^\star \|$ $(j\in\mathbb{N})$ ensure that, 
for all $j\in\mathbb{N}$,
\begin{align*}
&\left(1-\alpha_{n_j}\right) \sum_{i\in\mathcal{I}}  \left\|x_{n_j}^{(i-1)} - y_{n_j}^{(i)}\right\|^2
< 2 I N_1 \left(1-\alpha_{n_j} \right)\gamma_{n_j},\\
&\left(1-\alpha_{n_j}\right) \sum_{i\in\mathcal{I}}  \left\| y_{n_j}^{(i)}  - T^{(i)} \left(y_{n_j}^{(i)} \right) \right\|^2
< 2 I N_1 \left(1-\alpha_{n_j} \right)\gamma_{n_j},
\end{align*}
which, together with (C6) and (C7), implies that
$\lim_{j\to\infty} \|x_{n_j}^{(i-1)} - y_{n_j}^{(i)}\| = 0$ and
$\lim_{j\to\infty} \| y_{n_j}^{(i)}  - T^{(i)} (y_{n_j}^{(i)} ) \| = 0$ $(i\in \mathcal{I})$.

(ii) 
The same reasoning as in the proofs of Lemmas \ref{lem:3}(ii) and \ref{lem:4}(i) lead to $\lim_{j\to\infty} \|x_{n_j} - x_{n_j}^{(i-1)}\| = 0$, 
$\lim_{j\to\infty} \|x_{n_j} - y_{n_j}^{(i)}\| = 0$, and 
$\lim_{j\to\infty} \| x_{n_j}  - T^{(i)} (x_{n_j} ) \| = 0$ $(i\in \mathcal{I})$.
Assumption \ref{assum:KM-2} and $\lim_{j\to\infty} \|x_{n_j} - y_{n_j}^{(i)}\| = 0$ imply the boundedness of $(x_{n_j})_{j\in\mathbb{N}}$.

(iii)  
A discussion similar to the one for obtaining \eqref{4.2} means that, for all $j\in\mathbb{N}$,
$f (x_{n_j} ) - f(x_0^\star) < N_2 \sum_{i\in\mathcal{I}} \| x_{n_j} - y_{n_j}^{(i)} \|$, which, together with Lemma \ref{lem:4}(ii), means that
\begin{align*}
\limsup_{j\to \infty} \left(f \left(x_{n_j} \right) -  f \left(x_0^\star \right) \right) \leq N_2 \sum_{i\in\mathcal{I}} \lim_{j\to\infty} \left\| x_{n_j} - y_{n_j}^{(i)} \right\| = 0.
\end{align*}
Thus, $\limsup_{j\to \infty} f (x_{n_j} ) \leq f(x_0^\star)$.

(iv)
The boundedness of $(x_{n_j})_{j\in\mathbb{N}}$ implies that there exists $(x_{n_{j_k}})_{k\in\mathbb{N}} \subset (x_{n_j})_{j\in\mathbb{N}}$ such that $(x_{n_{j_k}})_{k\in\mathbb{N}}$ weakly converges to $x^\star$.
The same discussion as in the proof of Lemma \ref{lem:3-1}(ii) leads to $x^\star \in X$.
Moreover, (A2) and 
Proposition \ref{bc} imply that $f(x^\star) \leq \liminf_{k\to\infty} f(x_{n_{j_k}})$.
Accordingly, Lemma \ref{lem:4}(iii) guarantees that 
\begin{align*}
f(x^\star) \leq \liminf_{k\to\infty} f \left(x_{n_{j_k}}\right) \leq \limsup_{k\to\infty} f \left(x_{n_{j_k}}\right) 
\leq \limsup_{j\to \infty} f \left(x_{n_j} \right) \leq f\left(x_0^\star\right).
\end{align*}
That is, $x^\star \in X^\star$. 
From $\lim_{j\to\infty} \| x_{n_j} - x_{n_j}^{(i-1)} \| = 0$ $(i\in\mathcal{I})$, 
$(x_{n_{j_k}}^{(i)})_{k\in\mathbb{N}}$ $(i\in\mathcal{I})$ weakly converges to $x^\star \in X^\star$.
This completes the proof.
\qed

\section{Numerical Examples}\label{sec:5}
Consider the following problem with nonsmooth, convex objective functions \cite[Example 28]{comb2007} (see also Corollaries \ref{cor:1} and \ref{cor:2}).

\begin{prob}\label{concrete}
Assume that user $i$ ($i\in \mathcal{I}:= \{1,2,\ldots,I\}$) has its own private parameters $\omega_j^{(i)} > 0$, $a_j^{(i)} \in \mathbb{R}$,
$d_k^{(i)} \in \mathbb{R}$, and $c_k^{(i)} \in \mathbb{R}^N$ with $c_k^{(i)} \neq 0$, where $j\in \mathcal{N} := \{1,2,\ldots,N\}$ and
$k\in \mathcal{K} := \{1,2,\ldots,K\}$.
Define $f^{(i)} \colon \mathbb{R}^N \to \mathbb{R}$ and $C_k^{(i)} \subset \mathbb{R}^N$ ($i\in \mathcal{I}, k\in \mathcal{K}$) using
\begin{align*}
&f^{(i)}(x) := \sum_{j\in \mathcal{N}} \omega_j^{(i)} \left| x_j - a_j^{(i)} \right|  \text{ } \left(x\in \mathbb{R}^N \right) \text{ and } 
C_k^{(i)} := \left\{ x\in \mathbb{R}^N \colon \left\langle c_k^{(i)}, x  \right\rangle \leq d_k^{(i)}   \right\}.
\end{align*}
Then, 
\begin{align*}
\text{minimize } \sum_{i\in\mathcal{I}} f^{(i)} (x) 
\text{ subject to } x\in \bigcap_{i\in \mathcal{I}} C_{g^{(i)}},
\end{align*}
where 
$C_{g^{(i)}}$ ($i\in \mathcal{I}$)
is the generalized convex feasible set defined by \eqref{g} and \eqref{gcfs}
when $w_k^{(i)} := 1/K$ and $X^{(i)} = C := \{ x\in \mathbb{R}^N \colon \|x\| \leq 1 \}$ ($i\in \mathcal{I}, k\in \mathcal{K}$).  
\end{prob}

Here, $T^{(i)} \colon \mathbb{R}^N \to \mathbb{R}^N$ $(i\in\mathcal{I})$ is defined by \eqref{T} with $X^{(i)} = C$ and $w_k^{(i)} := 1/K$ $(k\in \mathcal{K})$.
Accordingly, $T^{(i)}$ $(i\in\mathcal{I})$ is firmly nonexpansive with 
$\mathrm{Fix}(T^{(i)}) = C_{g^{(i)}}$
(see section \ref{sec:3}). 
Hence, it is evident that Problem \ref{concrete} is an example of Problem \ref{problem:1}.

The experimental evaluations of the two proposed algorithms were done using a 27-inch iMac with a 3.2 GHz Intel Core i5 processor and 24 GB 1600 MHz DDR3 memory. 
The algorithms were written in Java 1.8.0\_60-b27
with $N := 100$, $I := 10$, and $K := 3$.
The values of $\omega_j^{(i)} \in (0,1]$, $a_j^{(i)} \in [-3,3]$, 
$d_k^{(i)} \in [0,1]$, $c_k^{(i)}$ with $\|c_k^{(i)}\| = 1$, and $x^{(i)}$ 
were randomly generated using org.apache.commons.math3.random.MersenneTwister. 
Algorithm \ref{algorithm:2} was used with 
\eqref{y_n_2} when $X^{(i)} := C$, 
and $(\alpha_n)_{n\in \mathbb{N}}$ and $(\gamma_n)_{n\in \mathbb{N}}$ were defined by\footnote{Numerical results in \cite{iiduka_siopt2013,iiduka_mp2015} indicate that the existing fixed point algorithms with small step sizes 
(e.g., $\gamma_n := 10^{-2}/(n+1)^a, 10^{-3}/(n+1)^a$) have faster convergence. Hence, the experiment described in this section used the step sizes in \eqref{sequence}.}
\begin{align}\label{sequence}
\gamma_n := \frac{10^{-3}}{(n+1)^a} \text{ } \left(a = \frac{1}{4}, \frac{1}{8}\right) \text{ and } 
\alpha_n := \frac{10^{-3}}{(n+1)^b} \text{ } \left( b = \frac{1}{2}, \frac{3}{4} \right)
\end{align}
while Algorithm \ref{algorithm:KM} was used with \eqref{y_n_2} when $x^{(i)}$ was replaced by $x_n^{(i-1)}$, $X^{(i)} := C$, $\alpha_n := t = 1/2$, and $(\gamma_n)_{n\in \mathbb{N}}$ was as given in \eqref{sequence}.

The incremental subgradient method (ISM) \cite{iiduka_2014} 
and parallel subgradient method (PSM) \cite{iiduka_fpta2015} 
were used for comparison.
ISM can be obtained by replacing $\mathrm{Prox}_{\gamma_n f^{(i)}}(x_n^{(i-1)})$ in Algorithm \ref{algorithm:KM} with $x_n^{(i-1)} - \gamma_n g_n^{(i)}$, where $g_n^{(i)} \in \partial f^{(i)}(x_n^{(i-1)})$.
The sequence generated by PSM is defined by 
$x_{n+1} := (1/I) \sum_{i\in \mathcal{I}} x_n^{(i)}$, where
$x_n^{(i)} := t x_n + (1-t) T^{(i)} (x_n - \gamma_n g_n^{(i)})$
and $g_n^{(i)}\in \partial f^{(i)}(x_n)$.
It is evident that Algorithms \ref{algorithm:2} and \ref{algorithm:KM} use the proximity operators of $f^{(i)}$s while ISM and PSM use the subgradients of $f^{(i)}$s.
To see how the choice of the order of the indices in $\mathcal{I}
:= \{1,2,\ldots,I\}$ 
affects the convergence rate of Algorithms \ref{algorithm:2} and \ref{algorithm:KM}, 
we compared Algorithms \ref{algorithm:2} and \ref{algorithm:KM} 
when (Case 1) $x_n^{(i)}$ $(i\in \mathcal{I})$ is calculated in the order of $1,2,\ldots,I$ and when (Case 2) $x_n^{(i)}$ $(i\in \mathcal{I})$ is calculated in randomly shuffled order. We found that the performances of Algorithms \ref{algorithm:2} and \ref{algorithm:KM} in Case 1 were almost the same as those in Case 2. Only the results for Case 1 are given due to lack of space.

One hundred samplings, each starting from a different randomly chosen initial point, were performed, and the results were averaged. 
Two performance measures were used. For each $n\in\mathbb{N}$,
\begin{align*}
F_n := \frac{1}{100} \sum_{s=1}^{100} \sum_{i\in\mathcal{I}} f^{(i)} \left(x_n (s) \right)
\text{ and }
D_n := \frac{1}{100} \sum_{s=1}^{100} \sum_{i\in\mathcal{I}} \left\| x_n(s) - T^{(i)} \left(x_n(s) \right)   \right\|, 
\end{align*}
where $(x_n(s))_{n\in\mathbb{N}}$ is the sequence generated from initial point $x(s)$ $(s=1,2,\ldots,100)$ for each of the four algorithms.
The value of $D_n$ represents the mean value of the sums of the distances between $x_n(s)$ and $T^{(i)}(x_n(s))$.
Hence, if $(D_n)_{n\in\mathbb{N}}$ converges to $0$, $(x_n)_{n\in\mathbb{N}}$ converges to some point in 
$\bigcap_{i\in\mathcal{I}} \mathrm{Fix}(T^{(i)}) = \bigcap_{i\in\mathcal{I}} C_{g^{(i)}}$.
\subsection{Case in which $C \cap \bigcap_{i\in \mathcal{I}} \bigcap_{k\in \mathcal{K}} C_k^{(i)} \neq \emptyset$}
Let us first consider Problem \ref{concrete} when 
the intersection of $C$ and $\bigcap_{i\in \mathcal{I}} \bigcap_{k\in \mathcal{K}} C_k^{(i)}$ is nonempty.

\begin{table}[H]
  \centering
  \caption{Comparison of proposed algorithms (Algorithms \ref{algorithm:2} and \ref{algorithm:KM}) with the existing algorithms (ISM and PSM) when $C \cap \bigcap_{i\in\mathcal{I}}\bigcap_{k\in \mathcal{K}} C_k^{(i)} \neq \emptyset$ (Algorithm(i) (resp. Algorithm(ii)) uses \eqref{sequence} with $a = 1/4$ and $b = 1/2$ (resp. $a = 1/8$ and $b = 3/4$))}
  \label{table1}
  \begin{tabular}{l||lll|lll}
  \hline
  & \multicolumn{3}{c|}{$|F_{n-1}-F_n|<10^{-3}$} & \multicolumn{3}{c}{$|D_{n-1}-D_n|<10^{-6}$}\\
  \hline
       & $n$ & time [s] & $F_n$ & $n$ & time [s] & $D_n$\\
  \hline
Alg.3.1(i) & 1850 & 0.028578 & 757.573942 & 1867 & 0.028840 & 0.003832 \\
Alg.3.1(ii) & 638 & 0.010049 & 749.866910 & 696 & 0.010949 & 0.003741 \\
Alg.4.1(i) & 1186 & 0.018333 & 749.866177 & 40 & 0.000665 & 0.000590 \\
Alg.4.1(ii) & 643 & 0.010307 & 749.641208 & 32 & 0.000555 & 0.001088 \\
ISM(i) & 1182 & 0.020897 & 749.870305 & 40 & 0.000781 & 0.000590 \\
ISM(ii) & 635 & 0.011474 & 749.649724 & 139 & 0.002594 & 0.002777 \\
PSM(i) & $\ge{}2000$ & 0.035836 & 761.479481 & 410 & 0.007589 & 0.000284 \\
PSM(ii) & $\ge{}2000$ & 0.035874 & 755.881753 & 331 & 0.006216 & 0.000697 \\
  \hline
  \end{tabular}
\end{table}

Table \ref{table1} shows the number of iterations $n$ and elapsed time when the algorithms (Algorithms \ref{algorithm:2} and \ref{algorithm:KM}, ISM, and PSM) satisfied $|F_{n-1} - F_n | < 10^{-3}$ and $|D_{n-1} - D_n| < 10^{-6}$. As shown, the $(F_n)_{n\in \mathbb{N}}$ generated by the incremental algorithms (Alg.\ref{algorithm:2}(ii), Alg.\ref{algorithm:KM}(ii), and ISM(ii)) using \eqref{sequence} with $a = 1/8$ and $b = 3/4$ converged faster than those (Alg.\ref{algorithm:2}(i), Alg.\ref{algorithm:KM}(i), and ISM(i)) using $a = 1/4$ and $b = 1/2$. Slowly diminishing step sizes such as $\gamma_n = 10^{-3}/(n+1)^{1/8}$ apparently affect the fast convergence of the algorithms. The number of iterations when PSM satisfied $|F_{n-1} - F_n|<10^{-3}$ was more than 2000, and PSM converged slowly compared with the incremental algorithms. The $(D_n)_{n\in \mathbb{N}}$ generated by all of the algorithms converged to $0$; i.e., the algorithms converged to a point in the constrained set in Problem \ref{concrete}. Alg.\ref{algorithm:KM}(i) and ISM(i) performed better than Alg.\ref{algorithm:2}(i) and PSM(i), and Alg.\ref{algorithm:2}(ii), Alg.\ref{algorithm:KM}(ii), and ISM(ii) had almost the same performance and converged faster than PSM(ii).
%
%
%
%
%
\subsection{Case in which $C \cap \bigcap_{i\in\mathcal{I}}\bigcap_{k\in \mathcal{K}} C_k^{(i)} = \emptyset$}
Next, let us consider Problem \ref{concrete} when the intersection of $C$ 
and $\bigcap_{i\in\mathcal{I}}\bigcap_{k\in \mathcal{K}} C_k^{(i)}$ is empty. 
Here, we assume that all users have the same $T^{(i)}$ to satisfy 
$\bigcap_{i\in \mathcal{I}} C_{g^{(i)}} =  \bigcap_{i\in \mathcal{I}} \mathrm{Fix}(T^{(i)}) \neq \emptyset$.
Accordingly, we consider the problem of minimizing $\sum_{i\in\mathcal{I}}f^{(i)}$ over $C_{g^{(i)}} \neq \emptyset$, where $C \cap \bigcap_{k\in \mathcal{K}} C_k^{(i)} = \emptyset$.

\begin{table}[H]
  \centering
  \caption{Comparison of proposed algorithms (Algorithms \ref{algorithm:2} and \ref{algorithm:KM}) with the existing algorithms (ISM and PSM) when $C \cap \bigcap_{k\in \mathcal{K}} C_k^{(i)} = \emptyset$ (Elapsed time of computing $x_{1999}$ in PSM(i) (resp. PSM(ii)) was 0.037320 [s] (resp. 0.036403 [s]), and $F_{1999}$ in PSM(i) (resp. PSM(ii)) was $977.37171882$ (resp. $955.43829899$))}
  \label{table3}
  \begin{tabular}{l||lll|lll}
  \hline
  & \multicolumn{3}{c|}{$|F_{n-1}-F_n|<10^{-3}$} & \multicolumn{3}{c}{$|D_{n-1}-D_n|<10^{-6}$}\\
  \hline
  & $n$ & time [s] & $F_n$ & $n$ & time [s] & $D_n$\\
  \hline
Alg.3.1(i) & 1219 & 0.020320 & 847.919647 & 250 & 0.004305 & 0.001109 \\
Alg.3.1(ii) & 1419 & 0.023099 & 737.425991 & 82 & 0.001473 & 0.001106 \\ 
Alg.4.1(i) & 1347 & 0.021702 & 879.651936 & 67 & 0.001217 & 0.000439 \\
Alg.4.1(ii) & 1807 & 0.029254 & 776.443722 & 36 & 0.000708 & 0.000795 \\
ISM(i) & 1348 & 0.024998 & 879.602066 & 43 & 0.001264 & 0.000485 \\
ISM(ii) & 1803 & 0.032661 & 776.644340 & 29 & 0.000707 & 0.000813 \\
PSM(i) & 6 & 0.000250 & 996.357860 & 70 & 0.001550 & 0.000288 \\
PSM(ii) & 6 & 0.000238 & 996.326212 & 65 & 0.001376 & 0.000491 \\
  \hline
  \end{tabular}
\end{table}

Table \ref{table3} shows the results for Algorithms \ref{algorithm:2} and \ref{algorithm:KM}, ISM, and PSM. Although Alg.\ref{algorithm:2}(ii), Alg.\ref{algorithm:KM}(ii), and ISM(ii) needed more iterations to satisfy $|F_{n-1}-F_n|<10^{-3}$ than Alg.\ref{algorithm:2}(i), Alg.\ref{algorithm:KM}(i), and ISM(i), Alg.\ref{algorithm:2}(ii), Alg.\ref{algorithm:KM}(ii), and ISM(ii) better optimized $\sum_{i\in \mathcal{I}}f^{(i)}$ than Alg.\ref{algorithm:2}(i), Alg.\ref{algorithm:KM}(i), and ISM(i). PSM converged slowly compared with the incremental algorithms, as also seen in Table \ref{table1}. All the algorithms converged to a point in $C_{g^{(i)}}$ in the early stages and, in particular, Alg.\ref{algorithm:2}(ii) ($F_{1419} \approx 737$), which is based on the Halpern fixed point algorithm, performed better than the algorithms based on the Krasnosel'ski\u\i-Mann fixed point algorithm. This is because the Halpern fixed point algorithm can minimize a certain convex function over the fixed point set of a nonexpansive mapping while the Krasnosel'ski\u\i-Mann fixed point algorithm can only find a fixed point. Since Problem \ref{concrete} is to minimize a convex function over the fixed point set of a nonexpansive mapping, Alg.\ref{algorithm:2}(ii) based on the Halpern algorithm is better suited for Problem \ref{concrete} than the algorithms based on the Krasnosel'ski\u\i-Mann algorithm.
%
%

\section{Conclusion and future work}\label{sec:6}
The problem of minimizing the sum of all users' nonsmooth, convex objective functions over the intersection of all users' fixed point sets in a Hilbert space was discussed, and two incremental proximal point algorithms were presented for solving the problem. One combines an incremental subgradient method with the Halpern fixed point algorithm, and the other is based on the Krasnosel'ski\u\i-Mann fixed point algorithm. Convergence analysis showed that, under certain assumptions, any weak sequential cluster point of the sequence generated by the Halpern-type algorithm is guaranteed to belong to the solution set of the problem and that there exists a weak sequential cluster point of the sequence generated by the Krasnosel'ski\u\i-Mann-type algorithm, which also belongs to the solution set. Numerical evaluations using concrete, nonsmooth, convex optimization problems showed the efficiency of the two algorithms.

Although nonsmooth, convex optimization with fixed point constraints in a Hilbert space was discussed, the numerically tested problems were defined in a finite-dimensional space. Future work includes generating numerical results that have special features of an infinite-dimensional space.

Since the bundle method \cite[chapter XIV]{hiri} is one of the most efficient methods for solving the problem of minimizing a general nonsmooth function, it would be of great interest to investigate whether bundle-type algorithms are well suited for nonsmooth (nonconvex) optimization with fixed point constraints. The first step would be to devise bundle-type algorithms for nonsmooth convex optimization over fixed point sets on the basis of previously reported results for the bundle method.

\textbf{Acknowledgments}
I am sincerely grateful to the editor, Immanuel Bomze, and the three anonymous reviewers for helping me improve the original manuscript. 
I also thank Kazuhiro Hishinuma for his input on the numerical examples.



\end{document}